\documentclass[preprint,10pt]{article}
\usepackage{amssymb}
\usepackage{mathrsfs}
\usepackage{amsfonts}
\usepackage{graphicx}
\usepackage{colortbl,dcolumn}
\usepackage{tabularx}
\usepackage{amsmath}
\usepackage{psfrag}
\usepackage{booktabs}
\usepackage{array}
\usepackage{enumerate}
\usepackage{cite}
\numberwithin{equation}{section}
\usepackage[top=1.2in, bottom=1.2in, left=0.9in, right=0.9in]{geometry}

\newcommand{\R}{\mathbb{R}}

\renewcommand{\P}{\mathbb{P}}

\newcommand{\dd}{\text{d}}
\newtheorem{theorem}{Theorem}[section]
\newtheorem{remark}[theorem]{Remark}
\newtheorem{assumption}[theorem]{Assumption}

\newtheorem{lemma}[theorem]{Lemma}

\newtheorem{Example}[theorem]{Example}

\begin{document}

%
 \title {Strong convergence rates of Galerkin finite element methods for  SWEs with cubic polynomial nonlinearity
\footnote{The authors would like to thank Dr. Jianbo Cui at HKPU for his helpful comments and suggestions.
R.Q. was supported by NSF of China (11701073)
and  by the Research Fund for  Yancheng Teachers University under 204040025.  X.W. was supported by NSF of China (12471394, 12071488, 12371417).
%
}
}

\author{
Ruisheng Qi$\,^\text{a}$,  \quad Xiaojie Wang$\,^\text{b}$ \\
\footnotesize $\,^\text{a}$ School of Mathematics and Statistics, Yancheng Teachers University, Yancheng, China\\
\footnotesize qiruisheng123@126.com\\
\footnotesize $\,^\text{b}$ School of Mathematics and Statistics, HNP-LAMA, Central South University, Changsha, China\\
\footnotesize x.j.wang7@csu.edu.cn\; and \;x.j.wang7@gmail.com
}
\maketitle
\begin{abstract}\hspace*{\fill}\\
  \normalsize
In the present work, strong approximation errors are analyzed for both the spatial semi-discretization and the spatio-temporal fully discretization of stochastic wave equations (SWEs)  with cubic polynomial nonlinearities and additive noises. The fully discretization is achieved by the standard Galerkin finite element method
in space and a novel exponential time integrator combined with the averaged vector field approach.  
The newly proposed scheme is proved to exactly satisfy a  trace formula  based on an energy functional. Recovering the convergence rates of the scheme, however, meets essential difficulties, due to the lack of the global monotonicity condition.  To overcome this issue, we derive the exponential integrability property of the considered numerical approximations, by the energy functional. Armed with these properties, we obtain the strong convergence rates of the approximations in both spatial and temporal direction. Finally, numerical results are presented to verify the previously theoretical findings.

\end{abstract}
\section{Introduction}

An instance of wave motions, involving the movement of a string \cite{cabana1972barrier} or a strand of DNA \cite{Dalang2009A}, is a commonly observed physical phenomenon that can be mathematically described by wave equations.
Due to the uncertainty of the media, stochastic effects are important in more accurately modeling. 
A stochastic perturbation of the deterministic model thus gives rise to stochastic wave equations, which are extensively studied in the literature \cite{anton2016full,cao2007spectral,chow2002stochastic,chow2006stochastic,cohen2022numerical,Cohen2013A,cox2019weak,cuistrong,Cui2022-semi,Dalang2009A,David2015fully,hong2022energy,jacobe2021weak,jin2022convergence,kovacs2020weak,KovFinite,li2017galerkin,li2022finite,Naurois2016Lower,PeszatNonlinear,Wang2014-higher,Yang2015-full,cai2023strong,lei2023numerical}.
In this paper,
we are interested in the numerical approximation of the semi-linear stochastic wave equations (SWEs)
\begin{align}\label{eq:sswe}
\begin{split}
\left\{\begin{array}{ll}
\dd u(t) =v(t)\dd t,&(t,x)\in [0, T]\times \mathcal{D},
\\
\dd v(t)= -\Lambda u(t) \dd t
-f(u(t)) \dd t + \dd W(t),& (t,x)\in [0, T]\times \mathcal{D},
\\
u(0)=u_0,\;v(0)=v_0,&x\in \mathcal{D},
\end{array}\right.
\end{split}
\end{align}
in a real separable Hilbert space $H:=L^2(\mathcal{D})$ with inner product $\left< \cdot,\,\cdot \right>$ and the associated norm $\|\cdot\|=\left<\cdot,\,\cdot\right>^{\frac12}$. 
Here 
$
-\Lambda:=\sum_{i=1}^d\frac{\partial^2 u}{\partial x^2_i}
$
is assumed to be the
Laplacian with Dirichlet boundary conditions
on a bounded convex domain $\mathcal{D}\subset \mathbb{R}^d$, $d = 1, 2$, with polygonal boundary. The stochastic process {$\{W(t)\}_{t\geq 0}$ is an $H$-valued $Q$-Wiener process on a filtered probability space
$(\Omega, \mathcal{F}, \mathbf{P}, \{\mathcal{F}_t\}_{t\geq 0})$ with respect to the normal filtration $\{\mathcal{F}_t\}_{t\geq 0}$. The initial data $u_0$ and $v_0$ are
random variables and the nonlinearity
$f(u)$ is assumed to be polynomial, which will be specified later.

In the last two decades, many researchers have conducted theoretical and numerical analysis of stochastic wave equations with global Lipschitz coefficients, see \cite{anton2016full,David2015fully,DalangH,Naurois2016Lower,PeszatNonlinear,Wang2014-higher,Cohen2013A,cohen2022numerical,kovacs2020weak,feng2022higher,jin2022convergence,jacobe2021weak,cao2007spectral,cox2019weak,kovacs2013weak,li2017galerkin,lei2023numerical,hausenblas2010weak}, just to mention a few.
In contrast to the globally Lipschitz case, 
stochastic wave equations with non-globally Lipschitz coefficients are less understood, from both theoretical and numerical points of view. 
Indeed, the well-posedness of SWEs with general polynomial drift coefficients has been well established in the literature \cite{chow2002stochastic,chow2006stochastic}.
For one-dimensional SWEs with cubic nonlinearity and $Q$-regular additive space-time noise, the author of \cite{SchurzAnalysis} proposed some nonstandard partial-implicit difference methods, which are proved to preserve the energy functional. But no strong convergence analysis was given there.
Recently, the authors of \cite{cuistrong} constructed  a different exponentially integrable scheme for SWEs with a cubic polynomial nonlinearity and an additive noise. 
Moreover, strong convergence rates of the proposed schemes were revealed, based on
the spatial-temporal regularity and exponential integrability properties of both the exact and numerical solutions.

For SWEs with super-linear coefficients and multiplicative noises, we are also aware of some recent works.
For example, the authors of \cite{li2022finite}
analyzed a finite element fully discrete scheme for multiplicative noise driven SWEs, but the strong convergence rate was only obtained in sub-spaces, not the whole probability space. 
In \cite{Cui2022-semi}, the authors proposed and analyzed stochastic scalar auxiliary variable (SAV) schemes,
which give the modified energy evolution law.
A strong convergence rate was derived there under the global Lipschitz condition.
In light of the discrete gradient method and the Pad\'{e} approximation,  \cite{hong2022energy} proposed several fully-discrete schemes for semi-linear
SWEs driven by multiplicative noises,
which preserve the averaged energy evolution law in the additive noise case. However, on convergence analysis was done there.

 
%

%
It is worth mentioning that the deterministic wave equation satisfies the energy conservation law.
Different from the deterministic case,
the SWE \eqref{eq:sswe} is not energy-conservative any more but satisfies the stochastic
energy evolution law \eqref{eq:thm-energy-evolution-law},  which is also called a trace formula in the literature, i.e., 
a linear drift of the expected value of the energy of the problem.
When discretizing the stochastic wave equation \eqref{eq:sswe}, an interesting question arises as to whether the resulting numerical approximations can still preserve the stochastic energy evolution law \eqref{eq:thm-energy-evolution-law}.
For the linear SWEs (i.e., $f(u)\equiv 0$), Cohen et al. gave a positive answer to the question by introducing an exponential time integrator that is called a trigonometric method in \cite{Cohen2013A}.
However, for semi-linear SWEs,  the usual exponential time integrator fail to exactly preserve the  stochastic energy evolution law \eqref{eq:thm-energy-evolution-law} \cite[Theorem 9]{anton2016full}. 
Recently, the authors of \cite{cuistrong} exploited the averaged vector field (AVF) to design schemes that preserve the energy evolution law. Moreover, Cui et al. \cite{Cui2022-semi} introduced a SAV approach based scheme for semi-linear SWEs, which preserves the modified energy evolution law. However, it still remains unclear whether any exponential-type scheme can be constructed that exactly preserves the energy evolution law.

In this paper we aim to introduce a new energy evolution law preserving fully discrete method of exponential-type for SWEs \eqref{eq:sswe}, done by the Galerkin finite element method in space and an exponential time integrator for time discretization.
To be more precise, let $V_h$, $h\in (0,1)$  be a sequence of finite dimensional sub-spaces of $ H_0^1(\mathcal{D})$ and $\tau=T/N$ with $N\in \mathbb{N}^+$ be the time step-size. 
Then the proposed fully discrete Galerkin finite element method for \eqref{eq:sswe} is to find $V_h$-valued stochastic processes $U_h^{n+1}, V_h^{n+1}$ such that
 \begin{align}
U_h^{n+1}=C_h(\tau)U_h^n+\Lambda_h^{-\frac12}S_h(\tau)V_h^n
-
\Lambda_h^{-1}\big(I-C_h(\tau)\big)P_h\int_0^1f\big(U^n_h+\theta(U^{n+1}_h-U^n_h)\big)\,\dd \theta,\label{eq:U-formulation-introduction}
\\
V_h^{n+1}
=
-\Lambda_h^{\frac12}S_h(\tau)U_h^n
+
C_h(\tau)V^n_h-\Lambda_h^{-\frac12}S_h(\tau)P_h\int_0^1f\big(U^n_h+\theta(U^{n+1}_h-U^n_h)\big)\,\dd \theta+P_h\delta W^n,
\label{eq:V-formulation-introduction}
\end{align}
where $\delta W^n := W(t_{n+1})-W(t_n)$. 
The discrete Laplace operator $-\Lambda_h$ and $L^2$-projection $P_h$ are defined by \eqref{eq:definition-discrete-A} and \eqref{eq:defn-Ph}, respectively.
In addition,
$C_h(\tau) := \cos(\tau \Lambda_h^{\frac12})$ and $S_h(\tau) := \sin(\tau \Lambda_h^{\frac12})$ are the cosine and sine operators.

The newly proposed scheme \eqref{eq:U-formulation-introduction}-\eqref{eq:V-formulation-introduction} 
differs from the usual exponential type time integrator in the evaluation of $f$. Here we exploit 
the averaged vector field (AVF) approach to   
ensure that 
the proposed method exactly preserves the energy evolution law and satisfies the following discrete trace formula (see Lemma \ref{eq-lem:discrete-trace-formula-full}):
\begin{align}\label{into:discrete-trace-formula-full}
\mathbb{E}
\big[
J(U_h^n,V_h^n)
\big]
=
\mathbb{E}\big[J(U_h^0,V_h^0)\big]
+
\tfrac12 \mathrm{Tr}(P_hQP_h)t_n,
\end{align}
where $J$ is a Lyaponov function defined by \eqref{eq:definition-J}.
The trace formula \eqref{into:discrete-trace-formula-full} can be achieved  because one has 
%
\begin{align}\label{eq-introd:energy-preseve}
J(U_h^{n+1},\overline{V}_h^{n+1})
=
J(U^n_h,V_h^n),\; with\; \overline{V}_h^{n+1}=V_h^{n+1}-P_h\delta W^n, \quad
\text{a.s..}
\end{align}
This identity suggests that 
the scheme \eqref{eq:U-formulation-introduction}-\eqref{eq:V-formulation-introduction} without the Wiener increment $\delta W^n$ is an energy-preserving method for the deterministic wave equation, which indeed reduces into the time-stepping scheme proposed by the authors of \cite{10.1093/imanum/dry047}. 
Owing to \eqref{eq-introd:energy-preseve}, one can obtain the exponential integrability properties of 
the numerical approximation,
which is 
crucial to obtaining the strong convergence rate
of the considered scheme.
It is shown that the scheme \eqref{eq:U-formulation-introduction}-\eqref{eq:V-formulation-introduction} 
is uniformly convergent (see Theorem \ref{them:covergence-full}):
\begin{align}
\sup_{n\in\{1,2,\cdots, N\}}\Big\|\|u(t_n)-U_h^n\|
+
\|v(t_n)-V_h^n\|_{\dot{H}^{-1}}\Big\|_{L^p(\Omega;\mathbb{R})}
\leq
C(h^{\frac{r+1}{r+2}\min\{\gamma,r+2\}}+\tau),
\end{align}
where 
the parameter $\gamma\in [1,\infty)$ for $d=1$ and $\gamma\in [1,3]$ for $d=2$, satisfying $\|\Lambda^{\frac{\gamma-1}2}Q^{\frac12}\|_{\mathrm{HS}}<\infty$ quantifies the spatial correlation of the noise process (Assumption \ref{assum:eq-noise} ) and $r$ comes from Assumption \ref{eq:assumption-rh-ph}. 
In addition, we denote $\mathbb{H}^\beta= \dot{H}^\beta\times \dot{H}^{\beta-1}$, $\beta\in \mathbb{R}$ with $\dot{H}^\beta:=D(\Lambda^{\frac\beta2})$.
%
Moreover, if $\Lambda_hP_h=P_h\Lambda_h$, then
\begin{align}
\sup_{n\in\{1,2,\cdots, N\}}\Big\|\|u(t_n)-U_h^n\|
+
\|v(t_n)-V_h^n\|_{\dot{H}^{-1}}\Big\|_{L^p(\Omega;\mathbb{R})}
\leq
C(h^{\min\{\gamma,r+1\}}+\tau).
\end{align}
%
%

We highlight that, distinct from Allen-Cahn type stochastic partitial differential equations (SPDEs), the drift of the SWE, as an SPDE system in the product space $\mathbb{H}^\beta$, does not obey the one-sided Lipschitz condition (also called  the global monotonicity condition in the literature). This leads to essential difficulties in obtaining the strong convergence rates, which are highly non-trivial. To overcome them, we rely on exponential integrability properties of both the exact and numerical solutions and successfully recover the strong convergence rates (see the proof of Theorems \ref{them:convergence-rate-semi}, \ref{them:covergence-full}).
Also, we mention that the exponential type scheme usually achieves higher convergence rates than the time discretizations based on rational approximation to the exponential function such as the implicit Euler and Crank–Nicolson (see \cite{KovFinite,Cohen2013A,anton2016full,wangAn2015}).
Finally, it is worthwhile to point out that, even when reducing into the deterministic case,  the above convergence results improves relevant results in \cite{10.1093/imanum/dry047}, where
the error analysis was carried out in a globally Lipschitz setting (see \cite[Assumption 3.1]{10.1093/imanum/dry047}).

The rest of this paper is organized as follows. 
In the next section,   
some necessary preliminaries and the well-posedness of the SWEs under consideration are presented.
In Section 3,  we consider the Galerkin finite element semi-discretization of SWEs and carry out strong error analysis.
In Section 4, we introduce a novel fully-discrete finite element method for SWEs and then derive the strong convergence result. Finally, we present numerical experiments to support our theoretical analysis in the previous section.

\section{The considered stochastic wave equations}

For a separable $\mathbb{R}$-Hilbert space $\left(H,\left<\cdot,\cdot\right>,\|\cdot\|\right)$,
we denote by $\mathcal{L}(H)$ the Banach space of all linear bounded operators from $H$ into $H$.
Let $\mathcal{L}_2 (H)$ represent
the Hilbert space consisting of Hilbert-Schmidt operators from $H$ into $H$, equipped with the scalar product and the norm
\begin{align}
\left<T_1,T_2\right>_{\mathcal{L}_2(H)}
=
\sum_{i\in \mathbb{N}}\left<T_1\phi_i,T_2\phi_i\right>,
\;
\|T\|_{\mathcal{L}_2(H)}
=
\sum_{i\in \mathbb{N}}\|T\phi_i\|^2,
\end{align}
where $\{\phi_i\}_{i\in \mathbb{N}}$ is an orthonormal basis of $H$. This definition turns out to be independent of the choice of the orthonormal basis. 
It is well-known that, if $T\in \mathcal{L}_2(H)$ and $L\in \mathcal{L}_1(H)$, then $TL$, $LT\in \mathcal{L}_2(H)$ and
\begin{align}
\|TL\|_{\mathcal{L}_2(H)}
\leq
\|T\|_{\mathcal{L}_2(H)}\|L\|_{\mathcal{L}_1(H)}
,
\;
 \|LT\|_{\mathcal{L}_2(H)}
 \leq
 \|T\|_{\mathcal{L}_2(H)}\|L\|_{\mathcal{L}_1(H)}.
 \end{align}

\subsection{Settings}
%
We start by introducing a self-adjoint, positive definite, linear operator $\Lambda$ with its domain $D(\Lambda):= H^2(\mathcal{D})\cap H_0^1(\mathcal{D})$.
Next, denote the Hilbert space $\dot{H}^\alpha := D(\Lambda^{\frac \alpha2})$, $\alpha\in \mathbb{R}$
with the inner product
\[
\left<x,y\right>_{\dot{H}^\alpha}
:=
\left<\Lambda^{\frac \alpha 2}x, \Lambda^{\frac \alpha 2}y\right>
=
\sum_{j=1}^\infty\lambda_j^\alpha\left <x,e_j\right>\left<y,e_j\right>,
\quad
\forall x,y \in \dot{H}^\alpha,
\]
and the associated norm
 $\|\Lambda^{\frac \alpha 2}x\|:=\big(\sum_{j=1}^\infty\lambda_j^\alpha\left<x,e_j\right>^2\big)^{\frac12}$, where
$\{(\lambda_j,e_j)\}_{j=1}^\infty$ are the eigenpairs of $\Lambda$ with orthonormal engenvectors in $L^2(\mathcal{D})$.
It is known that $\dot{H}^0=L^2(\mathcal{D})$ and $\dot{H}^1=H_0^1(\mathcal{D})$. Now we reformulate the stochastic problem \eqref{eq:sswe}
in
the following abstract form
\begin{align}\label{eq:abstract-form}
\dd X(t)=AX(t)\dd t+\mathbb{F}\big(X(t)\big)\,\dd t+\mathbb{B}\,\dd W(t),\;t\in[0,T],
\end{align}
where
\[
X_0=\biggl[
\begin{array}{c}
u_0\\
v_0
\end{array}
\biggl]
,
\quad
A=\biggl[
\begin{array}{c c}0 &I\\
-\Lambda&0
\end{array}
\biggl]
,
\quad
\mathbb{F}\big(X(t)\big)
=\biggl[
\begin{array}{c}
0\\-f\big(u(t)\big)
\end{array}
\biggl]
\;\;and\;\;
\mathbb{B}
=
\biggl[
\begin{array}{c}
0\\
I
\end{array}
\biggl].
\]
Here and below, by $I$ we denote the identity
operator.
 In addition, we define the product Hilbert space $\mathbb{H}^\alpha:=\dot{H}^\alpha\times \dot{H}^{\alpha-1}$, for $\alpha\in \mathbb{R}$, endowed with the inner product
\[\left<X,Y\right>_{\mathbb{H}^\alpha}
:=
\left<x_1,y_1\right>_{\dot{H}^\alpha}
+
\left<x_2,y_2\right>_{\dot{H}^{\alpha-1}}
,
\]
for any $X=\big[x_1,x_2\big]^\top,\;Y=\big[y_1,y_2\big]^\top$ and the norm $\|X\|_{\mathbb{H}^\alpha}:
=
\left<X,X\right>_{\mathbb{H}^\alpha}^{\frac12}.$
Then the operator $A$ generates a $C_0$-group $E(t)$, $t\in \R$ on the Hilbert space $\mathbb{H}^1$,
given by
\[
E(t)=
\biggl[
\begin{array}{c c}
C(t) &\Lambda^{-\frac12}S(t)
\\
-\Lambda^\frac12S(t)&C(t)
\end{array}
\biggl],
\]
where $C(t)=\cos(t \Lambda^{\frac12})$ and $S(t)=\sin(t\Lambda^\frac12)$ are the cosine and sine operators, respectively.
In the following lemma, we collect some results concerning
the temporal H\"{o}lder continuity of the sine and cosine operators, which can be found, for example, in \cite{anton2016full}.
\begin{lemma}\label{lem:eq-temporal-result-E}
For any $\alpha\in[0,1]$, there exists a constant $C=C(\alpha)$ such that
\begin{align}\label{eq:temporal-regularity-S(t)}
\|\big(S(t)-S(s)\big)\Lambda^{-\frac \alpha2}\|_{\mathcal{L}(H)}
\leq
C (t-s)^\alpha,\quad  \|\big(C(t)-C(s)\big)\Lambda^{-\frac \alpha2}\|_{\mathcal{L}(H)}
\leq
C (t-s)^\alpha,
\end{align}
and
\begin{align}\label{eq:temporal-regularity-C(t)}
\|\big(E(t)-E(s)\big)X\|_{\mathbb{H}^0}
\leq
C(t-s)^\alpha\|X\|_{\mathbb{H}^\alpha},
\end{align}
for all $t\geq s\geq 0$.
\end{lemma}

Next, we put some assumptions on the nonlinearity $F$, the noise process $W(t)$ and the initial data $X_0$.
\begin{assumption}\label{assum:Nonlinearity}
(Nonlinearity) Let $f$ be a cubic polynomial such that
\begin{align*}
\frac{\dd F(u)}{\dd u}=f(u)
\;\ and \;\;
a_1\|u\|_{L^4}^4-b_1\leq \int_\mathcal{D}F(u)\,\mathrm{d} x
\leq
 a_2\|u\|_{L^4}^4-b_2
\end{align*}
hold for some positive constants $a_1,a_2,b_1,b_2$ and  the following one-sided Lipschitz condition
\begin{align}\label{eq:one-side-L-condition}
-\left<f(u)-f(v),u-v\right>
\leq
c\|u-v\|^2,
\quad
\forall u,v \in H
\end{align}
holds for a positive constant $c>0$.
\end{assumption}
\begin{assumption}\label{assum:eq-noise}(Noise term) 
Let $\{W(t)\}_{t\in[0,T]}$ be  an $H$-valued $Q$-Wiener process on the stochastic basis $\big(\Omega,\mathcal{F},\mathbb{P},\{\mathcal{F}_t\}_{t\in[0,T]}\big)$, where    $Q\in \mathcal{L}(H)$ is a bounded, self-adjoint, positive semi-definite operator, satisfying
\begin{align}\label{eq-assum:noise}
\|\Lambda^{\frac{\gamma-1}2 }Q^{\frac12}\|_{\mathcal{L}_2(H)}<\infty,
\end{align}
for $\gamma\in [1,\infty)$ in the case $d=1$ or 
$\gamma\in [1,3]$ in the case $d=2$ .
\end{assumption}

\begin{assumption} \label{assum:initial-value-u0}
(Initial value). Let $X_0:\Omega\rightarrow \mathbb{H}^0$ be an $\mathcal{F}_0/\mathcal{B}(\mathbb{H}^0)$-measurable random variable and satisfy for a sufficiently
large number $p_0 \in \mathbb{N}$,
\begin{align}
\mathbf{E}[\|X_0\|_{\mathbb{H}^\gamma}^{p_0}]<\infty,
\end{align}
where the parameter $\gamma$ comes from \eqref{eq-assum:noise}.
\end{assumption}

\subsection{Regularity results of the mild solution}
This part is devoted to the well-posedness and spatio-temporal regularity  of the considered problem \eqref{eq:abstract-form}. 
The following existence, uniqueness and regularity results of  the mild solution to \eqref{eq:abstract-form} are quoted from \cite{cuistrong}.
\begin{theorem}\label{them:existen=expon-integ-result-mild}
Under Assumptions \ref{assum:Nonlinearity}-\ref{assum:initial-value-u0}, the underlying problem \eqref{eq:abstract-form} admits a unique mild solution given by
\begin{align}\label{eq:mild-solution}
X(t)
=
E(t)X_0
+
\int_0^tE(t-s)\mathbb{F}(X(s))\,\mathrm{d} s
+
\int_0^tE(t-s)\mathbb{B}\,\mathrm{d} W(s),\; t\in[0,T],
\end{align}
which satisfies the following estimate,  $\forall p\geq 1$
\begin{align}\label{eq:bound-u-H1}
\sup_{s\in[0,T]}\|u(s)\|_{L^p(\Omega;\dot{H}^1)}<\infty.
\end{align}
Moreover,
 there exits a constant $\alpha\geq \frac12 \mathrm{Tr}(Q)$ such that the mild solution \eqref{eq:mild-solution} enjoys the following exponential integrability property, for  any $c>0$
\begin{align}\label{lem-eq:expon-integ-property-mild}
\sup_{s\in[0,T]}\mathbb{E}\Big[\exp\Big(\frac{J(u(s),v(s))}{\exp(\alpha s)}\Big)\Big]
+
\mathbb{E}\Big[\exp\Big(\int_0^Tc\|u(s)\|_{L^6}^2\,\mathrm{d} s\Big)\Big]<\infty,
\end{align}
and
the energy evolution law
\begin{align}\label{eq:thm-energy-evolution-law}
\mathbb{E}\left[J(u(t),v(t))\right]
=
\mathbb{E}\left[J(u_0,v_0)\right]
+
\tfrac 12\mathrm{Tr}(Q)t,
\end{align}
where $J(u,v)$ is a Lyaponov function, given by
\begin{align}\label{eq:definition-J}
J(u,v)
=
\tfrac12\|\nabla u\|^2
+
\tfrac12\|v\|^2+F(u)
+
C_1,\; for\; C_1\geq b_1.
\end{align}
\end{theorem}

For two Hilbert
spaces $H^i, i = 1, 2$,  we introduce two operators $P_1$ and $P_2$ defined by
\begin{align}
P_ix=x_i,\;x=\big[x_1,x_2\big]^\top\in H_1\times H_2.
\end{align}
Equipped with Theorem \ref{them:existen=expon-integ-result-mild}, we can derive the following further regularity results  in the next theorem, in which some results were  shown in \cite{cuistrong}.
\begin{theorem}\label{lem:spatial-temporal-mild-solution}
Under Assumptions \ref{assum:Nonlinearity}-\ref{assum:initial-value-u0}, the mild solution of \eqref{eq:abstract-form} enjoys the following regularity,
\begin{align}\label{eq-them:spatial-regularity-mild}
\sup_{s\in[0,T]}\|u(s)\|_{L^p(\Omega;\dot{H}^\gamma)}<\infty,\forall p\geq 1,
\end{align}
and
\begin{align}\label{them-eq:temporal-regularity-u}
\|u(t)-u(s)\|_{L^p(\Omega;H)}
\leq
C|t-s|.
\end{align}
\end{theorem}
To prove Theorem \ref{lem:spatial-temporal-mild-solution}, we introduce some basic inequalities.
Recall first the following embedding inequalities,
\begin{align}\label{eq:embedding-equatlity-I}
\dot{H}^1\subset L^p(\mathcal{D}),\forall p\geq 1\quad \text{and} \quad \dot{H}^\delta \subset C(\mathcal{D};\mathbb{R}),
\quad
\text{ for } \delta>\tfrac d2, \;
d\in\{1,2\}.
\end{align}
With \eqref{eq:embedding-equatlity-I} at hand, one can further show, 
for any $x\in L^{\frac65}(\mathcal{D})$,
\begin{align}\label{eq:embedding-equatlity-IIII}
\begin{split}
\|\Lambda^{-\frac12} x\|
=
\sup_{v\in \dot{H}}\frac{\big|\big<x,\Lambda^{-\frac 12}v\big>\big|}{\|v\|}
\leq
\sup_{v\in \dot{H}}\frac{\|x\|_{L^{\frac65}} \|\Lambda^{-\frac 12}v\|_{L^6}}{\|v\|}
\leq
C\sup_{v\in \dot{H}}\frac{\|x\|_{L^{\frac65}} \|v\|}{\|v\|}
\leq
C\|x\|_{L^{\frac65}}.
\end{split}
\end{align}
Additionally, from \eqref{eq:embedding-equatlity-I} and the equivalence of the norms in $\dot{H}^\beta$ and $W^{\beta,2}$, for $\beta\in(0,1)$,   it follows that  $\delta>\frac d2$, $d=\{1,2\}$
\begin{align}\label{eq:bound-f-beta}
\begin{split}
C\|f(u)\|_{\dot{H}^\beta}
\leq
&
\|f(u)\|_{W^{\beta,2}}
\leq
C\|f(u)\|
+
C\left(\int_\mathcal{D}\int_\mathcal{D}\frac{|f(u(x))-f(u(y))|^2}{\|x-y\|^{2\beta+d}}\,\dd y\,\dd x\right)^{\frac12}
\\
\leq
&
C(1+\|u\|)(1+\|u\|_{C(\mathcal{D},\mathbb{R})}^2)
+
C(1+\|u\|_{C(\mathcal{D},\mathbb{R})}^2)\left(\int_\mathcal{D}\int_\mathcal{D}\frac{|u(x)-u(y)|^2}{\|x-y\|^{2\beta+d}}\,\dd y\,\dd x\right)^{\frac12}
\\
\leq
&
C(1+\|u\|+\|u\|_{W^{\beta,2}})(1+\|u\|_{C(\mathcal{D},\mathbb{R})}^2)
\\
\leq
&
C(1+\|u\|_{\dot{H}^\beta})(1+\|u\|_{\dot{H}^\delta}^2),\;\forall u\in \dot{H}^\delta.
\end{split}
\end{align}
Since  the norm $\|\cdot\|_2$ is  equivalent on $\dot{H}^2$ to
the standard Sobolev norm $\|\cdot\|_{H^2(\mathcal{D})}$ and $H^2(\mathcal{D})$
is an algebra,  one can find a constant $C >0$ such that, for any $f,g\in \dot{H}^2$,
\begin{align}\label{eq:algebra-properties-Hs}
\|fg\|_{H^2(\mathcal{D})}
\leq
C\|f\|_{H^2(\mathcal{D})}\|g\|_{H^2(\mathcal{D})}
\leq
C\| f\|_{\dot{H}^2}\,\| g\|_{\dot{H}^2}.
\end{align}
{\it Proof of Theorem \ref{lem:spatial-temporal-mild-solution}.}
The   spatial-temporal regularity results \eqref{eq-them:spatial-regularity-mild}-\eqref{them-eq:temporal-regularity-u}  were shown for $\gamma\geq1$ in  dimension one and for $\gamma\in\{1,2\}$ in  dimension two in  \cite{cuistrong}. Hence it remains to show \eqref{eq-them:spatial-regularity-mild} in dimension two  for $\gamma\in (1,2)\cup (2,3]$.
For the case $\gamma\in (1,2)$, by utilizing \eqref{eq:bound-f-beta} with $\beta=\frac{\gamma-1}2$ and $\delta=1+\frac{\gamma-1}3$ and using the boundness of the sine and cosine operators,  one can observe that
\begin{align}
\begin{split}
\|u(t)\|_{\dot{H}^\gamma}
\leq
&
\|P_1E(t)X_0\|_{\dot{H}^\gamma}+\int_0^t\|\Lambda^{\frac{\gamma}2}P_1E(t-s)\mathbb{F}(X(s))\|\,\dd s
+
\Big\|\int_0^t\Lambda^{\frac{\gamma}2}P_1E(t-s)\mathbb{B}\dd W(s)\Big\|
\\
\leq
&
\|X_0\|_{\mathbb{H}^\gamma}+\int_0^t\|\Lambda^{\frac{\gamma-1}2}f(u(s))\|\,\dd s
+
\Big\|\int_0^t\Lambda^{\frac{\gamma-1}2}S(t-s)\dd W(s)\Big\|
\\
\leq
&
\|X_0\|_{\mathbb{H}^\gamma}+
C\int_0^t\Big(1+\|\Lambda^{\frac{1+\frac{\gamma-1}3}2}u\|^2\Big)\|\Big(1+\Lambda^{\frac{\gamma-1}2}u\|\Big)\,\dd s
+
\Big\|\int_0^t\Lambda^{\frac{\gamma-1}2}S(t-s))\dd W(s)\Big\|
\\
\leq
&
\|X_0\|_{\mathbb{H}^\gamma}
+
C\int_0^t\Big(1+\|\Lambda^{\frac\gamma 2}u\|\|\Lambda^{\frac{4-\gamma}6}u\|\Big)\Big(1+\|\Lambda^{\frac{\gamma-1}2}u\|\Big)\,\dd s
+
\Big\|\int_0^t\Lambda^{\frac{\gamma-1}2}S(t-s)\dd W(s)\Big\|
\\
\leq
&
\|X_0\|_{\mathbb{H}^\gamma}
+
C\int_0^t\Big(1+\|u\|_{\dot{H}^{\gamma-1}}
+
\|u\|_{\dot{H}^\gamma}\|u\|_{\dot{H}^{\frac{4-\gamma}3}}
+
\|u\|_{\dot{H}^\gamma}\|u\|_{\dot{H}^{\frac{4-\gamma}3}}\|u\|_{\dot{H}^{\gamma-1}}\Big)\,\dd s
\\
&
+
\Big\|\int_0^t\Lambda^{\frac{\gamma-1}2}S(t-s)\dd W(s)\Big\|,
\end{split}
\end{align}
where in the fourth inequality we  also used the following Sobolev interpolation inequality
\begin{align}\label{eq:holder-inqualtiy}
\|\Lambda^{\frac{1+\frac{\gamma-1}3}2}u\|
\leq
\|\Lambda^{\frac{4-\gamma}6}u\|^{\frac12} \|\Lambda^{\frac\gamma2}u\|^{\frac12},
 \;for
  \;\gamma\in(1,2),\;u\in \dot{H}^\gamma.
\end{align}
Let  $\beta=\max\{\frac{4-\gamma}3,\gamma-1\}\in(0,1)$ and
\begin{align}
\Psi(t)
=
\|X_0\|_{\mathbb{H}^\gamma}
+
C\int_0^t(1+\|u\|_{\dot{H}^{\gamma-1}})\,\dd s
+
\Big\|\int_0^t\Lambda^{\frac{\gamma-1}2}S(t-s)\dd W(s)\Big\|.
\end{align}
Then,
\begin{align}
\|u(t)\|_{\dot{H}^\gamma}
\leq
\Psi(t)
+
C\int_0^t(\|u\|_{\dot{H}^\beta}+\|u\|^2_{\dot{H}^\beta})\|u\|_{\dot{H}^\gamma}\,\dd r.
\end{align}
Therefore,  we apply  Gronwall's inequality to obtain, for $\beta=\max\{\frac{4-\gamma}3,\gamma-1\}\in(0,1)$
\begin{align}
\begin{split}
\|u(t)\|_{\dot{H}^\gamma}
\leq
&
\Psi(t)
+
C\int_0^t \Psi(s)(\|u\|_{\dot{H}^\beta}+\|u\|_{\dot{H}^\beta}^2)
\exp\Big(C\int_s^t(\|u\|_{\dot{H}^\beta}+\|u\|_{\dot{H}^\beta}^2)\,\dd r\Big)\,\dd s
\\
\leq
&
\Psi(t)
+
C\int_0^t \Psi(s)(\|u\|_{\dot{H}^\beta}+\|u\|_{\dot{H}^\beta}^2)
\exp\Big(C\int_s^t\|u\|_{\dot{H}^\beta}^2\,\dd r\Big)\,\dd s.
\end{split}
\end{align}
Taking expectations and using \eqref{eq:bound-u-H1} as well as the H\"{o}lder inequality we get
\begin{align}\label{eq:u-H-gamma-12}
\begin{split}
\mathbb{E}[\|u(t)\|^p_{\dot{H}^\gamma}]
\leq
&
C\mathbb{E}\big(\Psi(t)\big)^p
+
C\int_0^t \mathbb{E} \left(\Psi(s)(\|u\|_{\dot{H}^\beta}+\|u\|_{\dot{H}^\beta}^2)\exp\Big(C\int_s^t\|u\|_{\dot{H}^\beta}^2\,\dd r\Big)\right)^p\,\dd s
\\
\leq
&
C\mathbb{E}\big(\Psi(t)\big)^p
+
C\int_0^t \left(\mathbb{E} \left(\Psi(s)\right)^{4p}(1+\mathbb{E}\|u\|^{8p}_{\dot{H}^1})\right)^{\frac14}
\left(\mathbb{E}\exp\Big(2pC\int_s^t\|u\|_{\dot{H}^\beta}^2\,\dd r\Big)\right)^{\frac12}\,\dd s
\\
\leq
&
C\mathbb{E}\big(\Psi(t)\big)^p
+
C\int_0^t \left(\mathbb{E} \left(\Psi(s)\right)^{4p}\right)^{\frac14}
\left(\mathbb{E}\exp\Big(2pC\int_s^t\|u\|_{\dot{H}^\beta}^2\,\dd r\Big)\right)^{\frac12}\,\dd s
.
\end{split}
\end{align}
Furthermore, using Jensen's inequality, Young's inequality,  \eqref{lem-eq:expon-integ-property-mild} and the fact that $\|\Lambda^{\frac\beta2}u\|\leq \|u\|^{1-\beta}\|\Lambda^{\frac12}u\|^{\beta}$, $\beta\in(0,1)$ helps show that for small $\epsilon>0$
\begin{align}
\begin{split}
\mathbb{E}\Big[\exp\Big(\int_s^t2pc\|\Lambda^{\frac\beta2}u\|^2\,\dd r\Big)\Big]
\leq
&
\sup_{t\in[0,T]}\mathbb{E}\Big[\exp\big(2pcT\|\Lambda^{\frac\beta2}u(t)\|^2\big)\Big]
\\
\leq
&
\sup_{t\in[0,T]}\mathbb{E}\Big[\exp\big(2pcT\|u\|^{2(1-\beta)}\|\Lambda^{\frac12}u\|^{2\beta)}\big)\Big]
\\
\leq
&
\sup_{t\in[0,T]}\mathbb{E}\Big[\exp\Big(\frac{\|u\|_{\dot{H}^1}^2}{2\exp(\alpha t)}\Big)
\exp\Big(\exp\Big(\frac{\beta }{1-\beta}\alpha T\Big)\|u\|^2(2pcT)^{\frac1{1-\beta}}2^{\frac \beta{1-\beta}}\Big)\Big]
\\
\leq
&
C(\epsilon)\sup_{t\in[0,T]}\mathbb{E}\Big[\exp\Big(\frac{\|u\|_1^2}{2\exp(\alpha t)}\Big)
\exp\big(\epsilon\|u\|^4_{L^4}\big)\Big]
\\
\leq
&
C\sup_{t\in[0,T]}\mathbb{E}\Big[\exp\Big(\frac{J(u,v)}{\exp(\alpha t)}\Big)\Big]
<
\infty.
\end{split}
\end{align}
In addition, utilizing \eqref{eq:bound-u-H1}, Assumption \ref{assum:initial-value-u0} and the Burkholder-Davis-Gundy inequality yields
\begin{align}
\begin{split}
\mathbb{E}\big(\Psi(t)\big)^p
\leq
&
C\mathbb{E}\|X_0\|_{\mathbb{H}^\gamma}^p
+
C\int_0^t(1+\mathbb{E}\|\Lambda^{\frac{\gamma-1}2}u\|^p)\,\dd s
+
C\mathbb{E}\Big\|\int_0^t\Lambda^{\frac{\gamma-1}2}S(t-s)\dd W(s)\Big\|^p
\\
\leq
&
CT+C\mathbb{E}\|X_0\|_{\mathbb{H}^\gamma}^p
+
C\int_0^t\mathbb{E}\|u\|_{\dot{H}^1}^p\,\dd s
+
C\left(\int_0^t\|\Lambda^{\frac{\gamma-1}2}S(t-s)Q^{\frac12}\|^2_{\mathcal{L}_2(H)}\dd s\right)^{\frac p2}
<\infty.
\end{split}
\end{align}
Substituting the above two estimates back into \eqref{eq:u-H-gamma-12} 
yields \eqref{eq-them:spatial-regularity-mild} for the case $\gamma\in(1,2)$.

 For the case $\gamma\in(2,3]$ in  dimension two, we use \eqref{eq:algebra-properties-Hs} and \eqref{eq-them:spatial-regularity-mild} for $\gamma=2$ to obtain
\begin{align}\label{eq:bound-f-gamma-1}
\sup_{s\in[0,T]}\|\Lambda^{\frac{\gamma-1}2}f(u(s))\|_{L^p(\Omega;H)}
\leq
\sup_{s\in[0,T]}\|\Lambda f(u(s))\|_{L^p(\Omega;H)}
\leq
C(1+\sup_{s\in[0,T]}\|u(s)\|^3_{L^{3p}(\Omega;\dot{H}^2)})<\infty,
\end{align}
  and then apply Assumption \ref{assum:initial-value-u0}, the boundedness of the sine and cosine operators and the Burkholder-Davis-Gundy inequality  to get
 \begin{align}
\begin{split}
\|u(t)\|_{L^p(\Omega;\dot{H}^\gamma)}
\leq
&
\|P_1E(t)X_0\|_{L^p(\Omega;\dot{H}^\gamma)}+\int_0^t\|P_1E(t-s)\mathbb{F}(X(s))\|_{L^p(\Omega;\dot{H}^\gamma)}\,\dd s
\\
&
+
\left\|\int_0^t P_1E(t-s)\mathbb{B}\,\dd W(s)\right\|_{L^p(\Omega;\dot{H}^\gamma)}
\\
\leq
&
\|X_0\|_{L^p(\Omega;\mathbb{H}^\gamma)}+
C\int_0^t\|\Lambda^{\frac{\gamma-1}2}f(u(s))\|_{L^p(\Omega;H)}\,\dd s
+
C\Big(\int_0^t\| \Lambda^{\frac{\gamma-1}2}S(t-s)Q^{\frac12}\|^2_{\mathcal{L}_2(H)}\,\dd s\Big)^{\frac12}
\\
\leq
&
C\|X_0\|_{L^p(\Omega;\mathbb{H}^\gamma)}
+
Ct\sup_{s\in[0,T]}\|\Lambda^{\frac{\gamma-1}2}f(u(s))\|_{L^p(\Omega;H)}
+
Ct^{\frac12}
\\
<
&\infty.
\end{split}
\end{align}
Now the proof of this theorem is 
accomplished. $\hfill\square$

As a direct consequence of Theorem \ref{lem:spatial-temporal-mild-solution}, the forthcoming lemma follows.
\begin{lemma} \label{lem:spatial-temporal-f}
Let Assumptions \ref{assum:Nonlinearity}-\ref{assum:initial-value-u0} be fulfilled. Then the following results hold
\begin{align}\label{eq-lem:spatial-regularity-f(u)}
\sup_{s\in[0,T]}
\|\Lambda^{\frac{\gamma-1}2}f\big(u(s)\big)\|_{L^p(\Omega;H)}
<\infty,
\end{align}
and
\begin{align}\label{eq-lem:temporal-regularity-f(u)}
\|\Lambda^{-\frac12}\big(f(u(t))-f(u(s))\big)\|_{L^p(\Omega;H)}
\leq
C|t-s|.
\end{align}
\end{lemma}

At the end of this section, we present a trace formula for the semi-linear stochastic wave equation \eqref{eq:sswe} (see also \cite[Proposition 5]{anton2016full} and \cite[Theorem  4.1]{cuistrong}).
\begin{lemma}
Let Assumptions \ref{assum:Nonlinearity}-\ref{assum:initial-value-u0} be fulfilled and the Hamiltonian function
$J(u,v)$ be given in Theorem \ref{them:existen=expon-integ-result-mild}. Then the mild solution of \eqref{eq:abstract-form} satisfies the trace formula
\begin{align}\label{eq:lemm-trace-formula-mild}
\mathbb{E}[J(u(t),v(t))]
=
\mathbb{E}[J(u(0),v(0))]
+
\tfrac12\mathrm{Tr}(Q) t.
\end{align}
\end{lemma}
\section{Galerkin finite element spatial semi-discretization}
In this section, we consider the spatial semi-discrete  approximation of the stochastic problem \eqref{eq:abstract-form} and show the  exponential integrability property of the semi-discrete problem, which is essentially used in the convergence analysis.
%

\subsection{Semi-discrete Galerkin finite element method}
%
This subsection is devoted to the semi-discrete  Galerkin finite element approximation of the stochastic problem \eqref{eq:abstract-form}. Before introducing  the semi-discrete  approximation, we present a brief overview of the relevant Galerkin finite element method.
 Let $V_h$, $h\in(0,1]$ be a sequence of finite dimensional subspaces of $ H_0^1(\mathcal{D})$.
The discrete Laplace operator $-\Lambda_h:V_h\rightarrow V_h$ is defined by
\begin{align}\label{eq:definition-discrete-A}
\left<\Lambda_hv_h,\chi_h\right>=a(v_h,\chi_h):=\big<\nabla v_h, \nabla \chi_h\big>,
\quad
\forall v_h,\;\chi_h\in V_h.
\end{align}
The above definition yields that $\Lambda_h$ is  self-adjoint and positive definite on $V_h$ and the following result holds
\begin{align}
\|\Lambda_h^{\frac12} v_h\|^2=a(v_h,v_h)=\|\Lambda^{\frac12}v_h\|^2=\|\nabla v_h\|^2.
\end{align}
In addition, we require two projection operators in order to discretize the stochastic problem \eqref{eq:sswe}. Firstly we define the Riesz projection operator $R_h:\dot{H}^1\rightarrow V_h$ as follows:
\begin{align}\label{eq:defn-Rh}
\left<\nabla R_h v,\,\nabla \chi_h\right>=\left<\nabla v,\,\nabla \chi_h\right>,
\quad
\forall v\in \dot{H}^1,
\:
\forall \chi_h\in V_h.
\end{align}
Secondly, the orthogonal projection operator $P_h: \dot{H}^{-1}\rightharpoonup V_h$ is defined by
\begin{align}\label{eq:defn-Ph}
\left<P_hv, \chi_h\right>=\left<v,\chi_h\right>,\;\forall\;v\in \dot{H}^{-1},\;\forall\;\chi_h\in V_h.
\end{align}
It is obvious that $\|P_hv\|\leq \|v\|,\;\forall v\in L^2(\mathcal{D})$.
Next we make an assumption on the error estimate for the Galerkin Riesz projection $R_h$ in the norm $\|\cdot\|$ and on the stability of the orthogonal projection $P_h$ in the norm $\|\cdot\|_1$, which
will be confirmed later by the usual finite element method and the spectral Galerkin method.
\begin{assumption}\label{eq:assumption-rh-ph}
For $V_h\subset H_0^1(\mathcal{D})$, there exists a constant $C$ such that
\begin{align}
\|R_hv-v\|\leq Ch^s\|v\|_{\dot{H}^s},\;for\;all\;v\in \dot{H}^s,\;s\in\{1,2,\cdots, r+1\},h\in(0,1), \;for\;some\;r\in \mathbb{N}^+,
\end{align}
and
\begin{align}
\|P_hv\|_{\dot{H}^1}
\leq
C\|v\|_{\dot{H}^1},\;for\;all\;v\in \dot{H}^1.
\end{align}
\end{assumption}

\begin{Example}(Finite element method) Let $\mathcal{T}_h$ a family of regular
and quasi-uniform meshes of the domain $\mathcal{D}$, where $h$ is the maximal mesh size.
Let $V_h$  be the space of all continuous functions, which are piecewise linear on $\mathcal{T}_h$ and zero on the boundary $\partial \mathcal{D}$. Then $V_h\subset H_0^1(\mathcal{D})$ and Assumption \ref{eq:assumption-rh-ph} holds for $r=1$, (see\cite[Lemma 1.1]{yan2005galerkin}).
\end{Example}

\begin{Example}(Spectral Galerkin method) Let $\{(\lambda_j,e_j)\}_{j=1}^\infty$ be the eigenpairs of the operator $\Lambda$. By setting $h:=\lambda^{-\frac12}_{J+1}$ for some $J\in \mathbb{N}^+$, the spaces $V_h$ are then defined as
\[
V_h:=span\{e_j,j=1,2,\cdots, J\}.
\]
In this case, Assumption \ref{eq:assumption-rh-ph}  holds for every $r\in \mathbb{N}^+$ (see \cite[Example 3.4]{kruse2014optimal}). In addition, the operators $\Lambda_h$ and $P_h$ can commute with each other.
\end{Example}

Under Assumption \ref{eq:assumption-rh-ph}, the operators $\Lambda$ and $\Lambda_h$ obey ( \cite[Theorem 4.4]{kovacs2012weak})
\begin{align}\label{eq:relation-A-Ah}
C_1\|\Lambda_h^{\frac \beta2}P_hv\|
\leq
\|\Lambda^{\frac \beta2}v\|
\leq
C_2
\|\Lambda_h^{\frac \beta2}P_hv\|,\;v\in \dot{H}^\beta,\;\beta\in[-1,1].
\end{align}

The semi-discrete Galerkin finite element method is to find $X_h(t)=\big[u_h(t),\;v_h(t)\big]^\top\in V_h\times V_h$ such that
\begin{align}\label{eq:semidsicrete}
\left\{\begin{array}{ll}
\dd X_h(t)=A_hX_h(t)\,\dd t+\mathbb{F}_h\big(X_h(t)\big)\,\dd t+\mathbb{B}_h\,\dd W(t),&\;t\in[0,T],
\\
X_h(0)=X_{0,h},&
\end{array}\right.
\end{align}
where
\[A_h:=\biggl[\begin{array}{c
c}0&I\\-\Lambda_h&0\end{array}\biggr],
\quad
 \mathbb{F}_h\big(X_h(t)\big):=
\biggl[
      \begin{array}{c}0 \\ -P_hf\big(u_h(t)\big)\end{array}
\biggr],
\;\;
\mathbb{B}_h:=\biggl[\begin{array}{c}0\\P_h\end{array}\biggr]
\;\: \text{ and }
X_{0,h}:
=
\biggl[\begin{array}{c}P_hu_0\\P_hv_0\end{array}\biggr].
\]
Here $A_h$ generates an analytic semigroup on the Hilbert space $V_h\times V_h$, give by
\[
E_h(t)=
\biggl[
\begin{array}{c c}
C_h(t) &\Lambda_h^{-\frac12}S_h(t)
\\
-\Lambda_h^\frac12S_h(t)&C_h(t)
\end{array}
\biggl],
\]
where $C_h(t)= \cos(t \Lambda_h^{\frac12})$ and $S_h(t)= \sin(t\Lambda_h^\frac12)$ are the cosine and sine operators, respectively.

Similarly to the continuous  case,  the semi-discrete problem \eqref{eq:semidsicrete} admits one unique mild solution, given by
\begin{align}\label{eq:mild-solution-semi}
X_h(t)=E_h(t)P_hX_0
+
\int_0^tE_h(t-s)\mathbb{F}_h\big(X_h(s)\big)\,\dd s
+
\int_0^tE_h(t-s)\mathbb{B}_h\,\dd W(s).
\end{align}


By the same arguments  of the proof of \cite[Corollary 3.1]{cuistrong},  we can show the exponential integrability property
of the solution to semi-discrete  problem \eqref{eq:mild-solution-semi}.
\begin{lemma}\label{lem:eq-spatial-temporal-regularity}
Suppose Assumptions \ref{assum:Nonlinearity}-\ref{assum:initial-value-u0} hold, then the following results  hold, for $p\geq 1$
and any constant $c>0$
\begin{align}
\label{lem:eq-temporal-regularity-semi}
\mathbb{E}
\bigg[
\sup_{s\in[0,T]}
\Big( J(u_h(s),v_h(s)) \Big) ^p
\bigg]
<
\infty,
\end{align}
and
\begin{align}\label{lem:eq-exponential-integ-properties-semi}
\mathbb{E}
\bigg[
\exp
\bigg(
\int_0^Tc\|u_h(s)\|_{L^6}^2\,\mathrm{d} s
\bigg)
\bigg]
<
\infty.
\end{align}
\end{lemma}
Its proof can be found in the Appendix.

Similar to the continuous problem \eqref{eq:abstract-form}, the exact solution of \eqref{eq:semidsicrete} also satisfies the following trace formula (cf. \cite[
Proposition 8]{anton2016full}).
\begin{lemma}
Let Assumptions \ref{assum:Nonlinearity}-\ref{assum:initial-value-u0} be fulfilled and $X_h(t)=[u_h(t),v_h(t)]^\top$ be the mild solution of the problem \eqref{eq:semidsicrete}. Then $X_h(t)$ satisfies the trace formula
\begin{align}
\mathbb{E}\big[J\big(u_h(t),v_h(t)\big)\big]=\mathbb{E}\big[J\big(u_h(0),v_h(0)\big)\big]
+
\tfrac 12\mathrm{Tr}\big(P_hQP_h\big)t.
\end{align}
\end{lemma}
\subsection{Strong convergence rates of the spatial-discretization}
The target of this part is to derive error estimates of the semi-discrete Galerkin finite element  method for the stochastic problem \eqref{eq:abstract-form}. The error analysis heavily relies on  the exponential integrability property
of the finite element approximation $u_h( \cdot )$ obtained in previous subsection and error estimates of the corresponding deterministic error operators as presented below.
\begin{lemma}
Denote $X_0=[u_0,v_0]^\top$ and let
\begin{align}
F_h(t)X_0:=\big(E(t)-E_h(t)P_h\big)X_0,
\end{align}
and
\begin{align}
K_h(t)X_0
:=\big(\Lambda^{\frac12}S(t)-\Lambda_h^{\frac12}S_h(t)P_h\big)u_0
+
\big(C(t)-C_h(t)P_h\big)v_0.
\end{align}
Under Assumption \ref{eq:assumption-rh-ph}, there exists a constant $C$ independent of the mesh size $h$ such that
\begin{align}\label{eq:error-semigroup}
\|F_h(t)X_0\|_{\mathbb{H}^0}
\leq
Ch^{\frac{{r+1}}{r+2}\beta}\|X_0\|_{\mathbb{H}^\beta},\;\beta\in[0,r+2],\;t\geq 0,
\end{align}
and
\begin{align}\label{eq:determinisitic-error-velocity}
\|K_h(t)X_0\|
\leq
 Ch^{\frac{r+1}{r+2}(\beta-1)}\|X_0\|_{\mathbb{H}^\beta}, \;\beta\in [1,r+3],\;t\geq 0.
 \end{align}
Furthermore, if $\Lambda_h P_h=P_h\Lambda_h$, then
\begin{align}\label{eq:error-semigroup-ii}
\|F_h(t)X_0\|_{\mathbb{H}^0}
\leq
Ch^\beta\|X_0\|_{\mathbb{H}^\beta},\;\beta\in[0,r+1],\;t\geq 0,
\end{align}
and
\begin{align}\label{eq:determinisitic-error-velocity-ii}
\|K_h(t)X_0\|
\leq
 Ch^{\beta-1}\|X_0\|_{\mathbb{H}^\beta}, \;\beta\in [1,r+2],\;t\geq 0.
 \end{align}
\end{lemma}

Equipped with the above estimates, we are well-prepared to do the error analysis. The next theorem states one main result of this section on the strong convergence rates of the spatial semi-discretization.
\begin{theorem}\label{them:convergence-rate-semi}
 Let Assumptions \ref{assum:Nonlinearity}-\ref{assum:initial-value-u0} and \ref{eq:assumption-rh-ph} be fulfilled. Then it holds
\begin{align}\label{them-eq:error-semi-norm-H}
\|X(t)-X_h(t)\|_{L^p(\Omega;\mathbb{H}^0)}
\leq
C h^{\frac{r+1}{r+2}\min\{\gamma,r+2\}},\;\forall\; t\in(0,T].
\end{align}
Additionally, for the velocity  in  dimension one,  it holds
\begin{align}\label{them-eq:vecolcity-error-semi-norm-H}
\|v(t)-v_h(t)\|_{L^p(\Omega;H)}
\leq
Ch^{\frac{r+1}{r+2}\min\{\gamma-1,r+2\}}, \;\forall\;t\in(0,T].
\end{align}
Furthermore, if $\Lambda_hP_h=P_h\Lambda_h$, it holds,
\begin{align}\label{them-eq:error-semi-norm-H-ii}
\|X(t)-X_h(t)\|_{L^p(\Omega;\mathbb{H}^0)}
\leq
C h^{\min\{\gamma,r+1\}},\;\forall\;t\in(0,T],
\end{align}
and for the velocity in  dimension one
\begin{align}\label{them-eq:vecolcity-error-semi-norm-H-ii}
\|v(t)-v_h(t)\|_{L^p(\Omega;H)}
\leq
Ch^{\min\{\gamma-1,r+1\}},\;\forall\;t\in(0,T].
\end{align}
\end{theorem}
\begin{remark}
For the spectral Galerkin method, Assumption \ref{eq:assumption-rh-ph}  holds for every $r\in \mathbb{N}^+$ and $\Lambda_hP_h=P_h\Lambda_h$. Therefore, the errors of the spectral Galerkin method can be estimated as  follows,
\begin{align}\label{them-eq:error-semi-norm-H-ii-i}
\|X(t)-X_h(t)\|_{L^p(\Omega;\mathbb{H}^0)}
\leq
C h^{\gamma},\;with\;h=\lambda_{J+1}^{-\frac12},\;\forall\;t\in(0,T],
\end{align}
and for the velocity in  dimension one
\begin{align}\label{them-eq:vecolcity-error-semi-norm-H-ii-ii}
\|v(t)-v_h(t)\|_{L^p(\Omega;H)}
\leq
Ch^{\gamma-1},\;with\;h=\lambda_{J+1}^{-\frac12}\;\forall\;t\in(0,T].
\end{align}
This comment also applies to the error estimate (Theorem \ref{them:covergence-full}) below for the fully-discretization.
\end{remark}
{\it Proof of Theorem \ref{lem:eq-spatial-temporal-regularity}.} In what follows we only present the proof of \eqref{them-eq:error-semi-norm-H} and \eqref{them-eq:vecolcity-error-semi-norm-H}. By the same arguments of the proof of \eqref{them-eq:error-semi-norm-H} and \eqref{them-eq:vecolcity-error-semi-norm-H}, we can show \eqref{them-eq:error-semi-norm-H-ii} and \eqref{them-eq:vecolcity-error-semi-norm-H-ii}. Let us begin with the proof of \eqref{them-eq:error-semi-norm-H}.
Subtracting \eqref{eq:mild-solution-semi} from \eqref{eq:mild-solution}   gives
\begin{align}\label{eq:split-semi-error}
\begin{split}
X(t)-X_h(t)
=
&
\big(E(t)-E_h(t)P_h\big)X_0
\\
&-
\int_0^t\big(E(t-s)\mathbb{F}(X(s))-E_h(t-s)\mathbb{F}_h(X_h(s))\big)\,\dd s
\\
&+
\int_0^t\big(E(t-s)\mathbb{B}-E_h(t-s)\mathbb{B}_h\big)\,\dd W(s)
\\
=&
I_1+I_2+I_3.
\end{split}
\end{align}
Subsequently we will deal with the three terms $I_1$, $I_2$, $I_3$, separately. By  \eqref{eq:error-semigroup} with $\beta=\min\{\gamma,r+2\}$,
 we bound the term $I_1$ as follows:
\begin{align}\label{eq:estimation-I1}
\|I_1\|_{\mathbb{H}^0}
\leq
Ch^{\frac{r+1}{r+2}\min\{\gamma,r+2\}}\|X_0\|_{\mathbb{H}^{\min\{\gamma,r+2\}}}.
\end{align}
In order to treat the term $I_2$, we first use \eqref{eq:relation-A-Ah} and \eqref{eq:embedding-equatlity-IIII}  to obtain
\begin{align}\label{eq:relation-H(-1)-L78}
\begin{split}
\|\Lambda_h^{-\frac12}P_h w\|
\leq
C\|w\|_{L^{\frac 65}},\;
\forall w\in H.
\end{split}
\end{align}
By applying the above estimate and \eqref{eq:error-semigroup} with $\beta=\min\{\gamma,r+2\}$, one can infer that,
\begin{align}\label{eq:estimation-I2}
\begin{split}
\|I_2\|_{\mathbb{H}^0}
\leq
&
\int_0^t\big\|\big(E(t-s)-E_h(t-s)P_h\big)\mathbb{F}(X(s))\|_{\mathbb{H}^0}\,\dd s
\\
&+
\int_0^t\|E_h(t-s)P_h\big(\mathbb{F}(X(s))-\mathbb{F}(X_h(s))\big)\|_{\mathbb{H}^0}\,\dd s
\\
\leq
&
Ch^{\frac{r+1}{r+2}\min\{\gamma,r+2\}}\int_0^t\big\|\Lambda^{\frac{\min\{\gamma,r+2\}-1}2}f(u(s))\|\,\dd s
+
C\int_0^t\|\Lambda_h^{-\frac12}P_h\big(f(u(s))-f(u_h(s))\big)\|\,\dd s
\\
\leq
&
Ch^{\frac{r+1}{r+2}\min\{\gamma,r+2\}}\int_0^t\big\|\Lambda^{\frac{\gamma-1}2}f\big(u(s)\big)\|\,\dd s
+
C\int_0^t\|f\big(u(s)\big)-f\big(u_h(s)\big)\|_{L^{\frac65}}\,\dd s
\\
\leq
&
Ch^{\frac{r+1}{r+2}\min\{\gamma,r+2\}}\int_0^t\big\|\Lambda^{\frac{\gamma-1}2}f(u(s))\|\,\dd s
+
C\int_0^t\|u(s)-u_h(s)\|(1+\|u(s)\|_{L^6}^2+\|u_h(s)\|_{L^6}^2)\,\dd s.
\end{split}
\end{align}
Finally, by plugging \eqref{eq:estimation-I1} and \eqref{eq:estimation-I2} back into
\eqref{eq:split-semi-error},
we can conclude
\begin{align}
\begin{split}
\|X(t)-X_h(t)\|_{\mathbb{H}^0}
\leq
\Phi(t)
+
C\int_0^t\|u(s)-u_h(s)\|(1+\|u(s)\|_{L^6}^2+\|u_h(s)\|_{L^6}^2)\
\dd s,
\end{split}
\end{align}
where we denote
\[
\Phi(t)
:=
Ch^{\frac{r+1}{r+2}\min\{\gamma,r+2\}}\Big(\|X_0\|_{\mathbb{H}^\gamma}
+
\int_0^t\|\Lambda^{\frac{\gamma-1}2}f(u(s))\|\,\dd s\Big)
+
\Big\|\int_0^t\big(E(t-s)-E_h(t-s)P_h\big)\mathbb{B}\,\dd W(s)\Big\|_\mathbb{H}.
\]
This together with Gronwall's inequality  helps us to obtain
\begin{align}
\begin{split}
\|X(t)-X_h(t)\|_\mathbb{H}
\leq
&
\Phi(t)
+
C\int_0^t \Phi(s)(1+\|u(s)\|_{L^6}^2+\|u_h(s)\|^2_{L^6})
\\
&
\cdot
\exp\Big(c\int_s^t(1+\|u(r)\|_{L^6}^2+\|u_h(r)\|^2_{L^6})\,\dd r\Big)\,\dd s.
\end{split}
\end{align}
As a consequence,
\begin{align}
\begin{split}
\mathbb{E}\big[\|X(t)-X_h(t)\|^p_{\mathbb{H}^0}\big]
\leq
&
\mathbb{E}[|\Phi(t)|^p]
+
C\int_0^t\big(\mathbb{E}[|\Phi(s)|^{4p}])^{\frac14}\big(\mathbb{E}[(1+\|u(s)\|_{L^6}^2+\|u_h(s)\|^2_{L^6})^{4p}]\big)^{\frac14}
\\
&\cdot
\Big(\mathbb{E}\Big[\exp\Big(2pc\int_s^t(1+\|u(r)\|_{L^6}^2+\|u_h(r)\|^2_{L^6})\,\dd r\Big)\Big]\Big)^{\frac12}\,\dd s
\\
\leq
&
\mathbb{E}[|\Phi(t)|^p]
+
C\int_0^t\big(\mathbb{E}[|\Phi(s)|^{4p}])^{\frac14}\,\dd s,
\end{split}
\end{align}
where in the second inequality we also used  \eqref{lem-eq:expon-integ-property-mild} and \eqref{lem:eq-exponential-integ-properties-semi}.
Furthermore,  using \eqref{eq-lem:spatial-regularity-f(u)}, \eqref{eq:error-semigroup} with $\beta=\min\{\gamma,r+2\}$ and the
Burkholder-Davis-Gundy inequality leads to
\begin{align}
\begin{split}
\mathbb{E}[|\Phi(t)|^p]
\leq
&
Ch^{\frac{r+1}{r+2}p\min\{\gamma,r+2\}}\mathbb{E}[\|X_0\|^p_{\mathbb{H}^\gamma}]
+
Ch^{\frac{r+1}{r+2}p\min\{\gamma,r+2\}}\int_0^t\mathbb{E}\big[\|\Lambda^{\frac{\gamma-1}2}f(u(s))\|\big]^p\,\dd s
\\
&+
C\Big(\int_0^t\|(E(t-s)-E_h(t-s)P_h)\mathbb{B}Q^{\frac12}\|_{\mathcal{L}_2(\mathbb{H}^0)}^2\,\dd s\Big)^{\frac p2}
\\
\leq
&
Ch^{\frac{r+1}{r+2}p\min\{\gamma,r+2\}}\Big(\mathbb{E}[\|X_0\|^p_{\mathbb{H}^\gamma}]
+
\sup_{s\in[0,T]}\mathbb{E}\big[\|\Lambda^{\frac{\gamma-1}2}f(u(s))\|\big]^p
+
\|\Lambda^{\frac{\min\{\gamma,r+2\}-1}2}Q^{\frac12}\|_{\mathcal{L}_2(H)}^p\Big)
\\
\leq
&
Ch^{\frac{r+1}{r+2}p\min\{\gamma,r+2\}}.
\end{split}
\end{align}
Hence
\[
\mathbb{E}[\|X(t)-X_h(t)\|_{\mathbb{H}^0}^p]\leq
Ch^{\frac{r+1}{r+2}p\min\{\gamma,r+2\}},\]
which shows  \eqref{them-eq:error-semi-norm-H}.
Now we turn our attention to the error $\|v(t)-v_h(t)\|_{L^p(\Omega;H)}$  in  dimension one. By \eqref{eq:split-semi-error} and definition of $E(t)$ and $E_h(t)$, we have
\begin{align}\label{eq:split-error-v-semi}
\begin{split}
\|v(t)-v_h(t)\|_{L^p(\Omega;H)}
\leq
&
\|P_2\big(E(t)-E_h(t)P_h\big)X_0\|_{L^p(\Omega;H)}
\\
&
+
\int_0^t\big\|P_2\big(E(t-s)\mathbb{F}(X(s))-E_h(t-s)P_h\mathbb{F}(X_h(s))\big)\big\|_{L^p(\Omega;H)}
\,\dd s
\\
&
+
\Big\|\int_0^tP_2\big(E(t-s)-E_h(t-s)P_h\big)\mathbb{B}\,\dd W(s)\Big\|_{L^p(\Omega;H)}
\\
&=L_1+L_2+L_3.
\end{split}
\end{align}
Using \eqref{eq:determinisitic-error-velocity} with $\beta=\min\{\gamma,r+3\}$ gives
\begin{align}\label{eq:estimate-L1}
\begin{split}
L_1= \|P_2F_h(t)X_0\|_{L^p(\Omega;H)}
\leq
Ch^{\frac{r+1}{r+2}\min\{\gamma-1,r+2\}}\|X_0\|_{L^p(\Omega;\mathbb{H}^\gamma)}.
\end{split}
\end{align}
In order to deal with the terms $L_2$, we decompose it into two parts
\begin{align}
\begin{split}
L_2
\leq
&
\int_0^t\big\|P_2\big(E(t-s)-E_h(t-s)P_h)\mathbb{F}(X(s)\big)\big\|_{L^p(\Omega;H)}
\,\dd s
\\
&
+
\int_0^t\big\|P_2E_h(t-s)P_h\big(\mathbb{F}(X(s))-\mathbb{F}(X_h(s))\big)\big\|_{L^p(\Omega;H)}\,\dd s
\\
=
&L_{21}+L_{22}.
\end{split}
\end{align}
Owing to \eqref{eq:determinisitic-error-velocity} with $\beta=\min\{\gamma,r+3\}$ and \eqref{eq-lem:spatial-regularity-f(u)}, we infer
\begin{align}\label{eq:estimate-L21}
\begin{split}
L_{21}
\leq
&
Ch^{\frac{r+1}{r+2}\min\{\gamma-1,r+2\}}\int_0^t\big\|\Lambda^{\frac{\min\{\gamma,r+3\}-1}2}f(u(s))\big\|_{L^p(\Omega;H)}
\,\dd s
\\
\leq
&
Ch^{\frac{r+1}{r+2}\min\{\gamma-1,r+2\}}t\sup_{s\in[0,T]}\|\Lambda^{\frac{\gamma-1}2}f(u(s))\big\|_{L^p(\Omega;H)}
\\
\leq
&
Ch^{\frac{r+1}{r+2}\min\{\gamma-1,r+2\}}.
\end{split}
\end{align}
Likewise, we use
 \eqref{them-eq:error-semi-norm-H},
  \eqref{eq:embedding-equatlity-I}, \eqref{eq:bound-u-H1} and \eqref{lem:eq-temporal-regularity-semi}  to arrive at
\begin{align}\label{eq:estimate-L22}
\begin{split}
L_{22}
\leq
&
\int_0^t\|P_2E_h(t-s)P_h\big(\mathbb{F}(X(s))-\mathbb{F}(X_h(s))\big)\|_{L^p(\Omega;H)}\,\dd s
\\
\leq
&
\int_0^t\|f(u(s))-f(u_h(s))\|_{L^p(\Omega;H)}\,\dd s
\\
\leq
&
C\int_0^t\|u(s)-u_h(s)\|_{L^{2p}(\Omega;H)}(1+\|u(s)\|_{L^{4p}(\Omega;C(\mathcal{D},\mathbb{R}))}^2
+
\|u_h(s)\|_{L^{4p}(\Omega;C(\mathcal{D},\mathbb{R}))}^2)\,\dd s
\\
\leq
&
Ch^{\frac{r+1}{r+2}\min\{\gamma,r+2\}}\big(1
+
\sup_{s\in[0,T]}\|u(s)\|_{L^{4p}(\Omega;\dot{H}^1)}^2
+
\sup_{s\in[0,T]}\|u_h(s)\|_{L^{4p}(\Omega;\dot{H}^1)}^2\big)
\\
\leq
&Ch^{\frac{r+1}{r+2}\min\{\gamma,r+2\}},
\end{split}
\end{align}
which together with \eqref{eq:estimate-L21} implies
\begin{align}\label{eq:estimate-L2}
L_2\leq Ch^{\frac{r+1}{r+2}\min\{\gamma-1,r+2\}}.
\end{align}
At last, the
Burkholder-Davis-Gundy inequality and \eqref{eq:determinisitic-error-velocity} with $\beta=\min\{\gamma,r+3\}$ help us to obtain
\begin{align}\label{eq:estimate-L3}
\begin{split}
L_3
\leq
\Big(\int_0^t\|(C(t-s)-C_h(t-s)P_h)Q^{\frac12}\|^2_{\mathrm{HS}}\,\dd s\Big)^{\frac 12}
\leq
C h^{\frac{r+1}{r+2}\min\{\gamma-1,r+2\}} \|\Lambda^{\frac{\min\{\gamma,r+3\}-1}2}Q^{\frac12}\|_{\mathrm{HS}}.
\end{split}
\end{align}
Now gathering \eqref{eq:estimate-L1}, \eqref{eq:estimate-L2} and \eqref{eq:estimate-L3}   together gives the estimate of $\|v(t)-v_h(t)\|_{L^p(\Omega;H)}$. The proof is thus complete. $\hfill\square$
\begin{remark}
For the case $d=3$, the strong convergence rate can not be derived, since we can not obtain the exponential integrability in $L^6$ for the solutions of stochastic problem \eqref{eq:sswe} and the semidiscrete problem \eqref{eq:semidsicrete} in this case. For more details, please refer to \cite{cuistrong}. In addition, we can not prove the strong convergence rate of the velocity in dimension two due to the absence of the boundedness  of $u_h(t)$ given by \eqref{eq:semidsicrete} in $C(\mathcal{D},\mathbb{R})$. Also, this comment applies to the error estimates (Theorem \ref{them:covergence-full}) below for the fully discretization.
\end{remark}
\section{Galerkin finite element fully discretization}

In this section, we propose and analyze a fully-discrete Galerkin finite element  method for the nonlinear stochastic wave equations \eqref{eq:abstract-form}.

\subsection{Fully discrete  Galerkin finite element method}
In order to introduce the fully discrete scheme, we first recall the mild solution of \eqref{eq:mild-solution-semi}
and then obtain
\begin{align}\label{eq:relation-mild-semidiscrete}
X_h(t_{n+1})=E_h(\tau )X_h(t_n)
+
\int_{t_n}^{t_{n+1}}
E_h(t_{n+1}-s)\mathbb{F}_h(X_h(s))\,\dd s
+
\int_{t_n}^{t_{n+1}}E_h(t_{n+1}-s)\mathbb{B}_h\,\dd W(s),
\end{align}
where $\tau$ is the time-step size, $t_n=n\tau$, for $n\in \{1,2,\cdots, N\}, N\in \mathbb{N}^+$.
Replacing $f(u_h(s))$ with the well-known AVF method for a function $\nabla_{AVF}f(u_h(t_n),u_h(t_{n+1})):=\int_0^1f\big(u_h(t_n)+\theta(u_h(t_{n+1})-u_h(t_n))\big)\,\dd \theta$,
the deterministic integral appearing in \eqref{eq:relation-mild-semidiscrete} can be approximated by
\begin{align}
\begin{split}
\int_{t_n}^{t_{n+1}} E_h(t_{n+1}-s)P_h\mathbb{F}(X_h(s))\,\dd s
\thickapprox
&
\biggl[
      \begin{array}{c}-\int_{t_n}^{t_{n+1}}\Lambda_h^{-\frac12}S_h(t_{n+1}-s)P_h\nabla_{AVF}f\big(u_h(t_n),u_h(t_{n+1})\big)\,\dd s \\-\int_{t_n}^{t_{n+1}} C_h(t_{n+1}-s)P_h\nabla_{AVF}f\big(u_h(t_n),u_h(t_{n+1})\big)\,\dd s\end{array}
\biggr]
\\
=
&
\biggl[
      \begin{array}{c}-\Lambda_h^{-1}(1-C_h(\tau))P_h\nabla_{AVF}f\big(u_h(t_n),u_h(t_{n+1})\big)
      \\
      -\Lambda_h^{-\frac12}S_h(\tau)P_h\nabla_{AVF}f\big(u_h(t_n),u_h(t_{n+1})\big)\end{array}
\biggr].
\end{split}
\end{align}
Based on the above observation, we introduce the fully discrete Galerkin finite element method for the stochastic wave equation \eqref{eq:abstract-form} as follows
\begin{align}
U_h^{n+1}=C_h(\tau)U_h^n+\Lambda_h^{-\frac12}S_h(\tau)V_h^n
-
\Lambda_h^{-1}(I-C_h(\tau))P_h\int_0^1f\big(U^n_h+\theta(U^{n+1}_h-U^n_h)\big)\,\dd \theta,\label{eq:U-formulation}
\\
V_h^{n+1}
=
-\Lambda_h^{\frac12}S_h(\tau)U_h^n+C_h(\tau)V^n_h
-
\Lambda_h^{-\frac12}S_h(\tau)P_h\int_0^1f\big(U^n_h+\theta(U^{n+1}_h-U^n_h)\big)\,\dd \theta+P_h\delta W^n,
\label{eq:V-formulation}
\end{align}
where we denote $\delta W^n := W(t_{n+1})-W(t_n)$. Let 
$X_h^n:=[U_h^n,V_h^n]^\top$, $\mathbb{B}_h=[0,P_h]^\top$, and
\[
\mathbb{F}^n_h
:=
\big[0,\;
      -\int_0^1P_hf\big(U_h^n+\theta(U_h^{n+1}-U_h^n)\big)\,\dd \theta\big]^\top.
\]
Then, by the fact
\[\int_{t_n}^{t_{n+1}}\Lambda_h^{-\frac12}S_h(t_{n+1}-s)\,\dd s
=
\Lambda_h^{-1}(1-C_h(\tau))\; and \;
\int_{t_n}^{t_{n+1}}C_h(t_{n+1}-s)\,\dd s
=
\Lambda_h^{-\frac12}S_h(\tau),
\]
the scheme \eqref{eq:U-formulation}-\eqref{eq:V-formulation} also can be rewritten as
\begin{align}\label{eq:solution-full}
X_h^{n+1}=E_h(\tau)X_h^{n}+\int_{t_n}^{t_{n+1}}E_h(t_{n+1}-s)\mathbb{F}^{n}_h\,\dd s
+
\int_{t_n}^{t_{n+1}}E_h(0)\mathbb{B}_h\,\dd W(s).
\end{align}
%
\begin{lemma}\label{lem:existence-unique-solution}
Suppose Assumptions \ref{assum:Nonlinearity}-\ref{assum:initial-value-u0} hold and $\tau<2$, then the fully discrete  scheme \eqref{eq:U-formulation}-\eqref{eq:V-formulation} is uniquely solvable almost surely.
\end{lemma}
{\it Proof of Lemma \ref{lem:existence-unique-solution}.} To show the existence and uniqueness of the fully discrete problem \eqref{eq:U-formulation}-\eqref{eq:V-formulation}, it suffices to  consider the following  auxiliary  problem
\begin{align}\label{eq:auxiliary-porblem-III}
\begin{split}
 \Lambda_h\chi'
=
-2\Lambda_h\sin^2(\tfrac{\tau\Lambda_h^{\frac12}}2)U_h^n
+
2\sin(\tfrac{\tau\Lambda_h^{\frac12}}2)\cos(\tfrac{\tau\Lambda_h^{\frac12}}2)\Lambda_h^{\frac12}V_h^n
-
2\sin^2(\tfrac{\tau\Lambda_h^{\frac12}}2)P_h\int_0^1f(U_h^n+\chi'\theta )\,\dd \theta.
\end{split}
\end{align}
Let $V_h=V_{h,0}\oplus V_{h,0}^T$ and $\widetilde{P}_h$ be a projection from $H$ to $V_{h,0}^T$, where $V_{h,0}$ is the zero space of the operator $\sin(\tfrac{\tau\Lambda_h^{\frac12}}2)$ in $V_h$. Then $\sin(\tfrac{\tau\Lambda_h^{\frac12}}2)$ is a bijection operator from the space $V_{h,0}^T$ to $V_{h,0}^T$ and the problem  \eqref{eq:auxiliary-porblem-III} is equivalent to the following problem
\begin{align}\label{eq:auxiliary-porblem-IIII}
\begin{split}
     \Lambda_h\chi'
    =
-2\Lambda_h\sin^2(\tfrac{\tau\Lambda_h^{\frac12}}2)\widetilde{P}_hU_h^n
+
2\sin(\tfrac{\tau\Lambda_h^{\frac12}}2)\cos(\tfrac{\tau\Lambda_h^{\frac12}}2)\Lambda_h^{\frac12}\widetilde{P}_hV_h^n
-
2\sin^2(\tfrac{\tau\Lambda_h^{\frac12}}2)\widetilde{P}_h\int_0^1f(U_h^n+\chi'\theta )\,\dd \theta.
\end{split}
\end{align}
Hence, it suffices to solve the following problem in $V_{h,0}^T$
\begin{align}\label{eq:auxiliary-porblem-I}
\begin{split}
\Lambda_h\chi
=
-\Lambda_h\sin(\frac{\tau\Lambda_h^{\frac12}}2)\widetilde{P}_hU_h^n
+
\cos(\frac{\tau\Lambda_h^{\frac12}}2)\Lambda_h^{\frac12}\widetilde{P}_hV_h^n
-
\sin(\frac{\tau\Lambda_h^{\frac12}}2)\widetilde{P}_h\int_0^1f(U_h^n+2\sin(\frac{\tau\Lambda_h^{\frac12}}2)\chi\theta )\,\dd \theta.
\end{split}
\end{align}

Denote 
\begin{equation}
G_h(\chi)
=
\Lambda_h\chi-b+\sin(\frac{\tau\Lambda_h^{\frac12}}2)
\widetilde{P}_h\int_0^1f\big(c+2\sin(\frac{\tau\Lambda_h^{\frac12}}2) \theta \chi\big)\,\dd \theta
\end{equation}
for any $ b\in V_{h,0}^T, c\in V_h$. 
Then $G:V_{h,0}^T\rightarrow V_{h,0}^T$.
It is a well-known simple consequence of Brouwer's fixed point
theorem that the equation $G_h(\chi) = 0$ has a solution $\chi \in B_q = \{\chi\in V_{h,0}^T; \|\chi\|\leq q\}$ if $\big<G(\chi), \chi\big> >0$ for $\|\chi\|=q$.
 In fact, if we assume that
$G(\chi) \neq 0$  in $B_q$, then the mapping $\Phi(\chi) = -\frac{q G_h(\chi)}{\|G_h(\chi)\|}:B_q\rightarrow B_q$
is continuous, and therefore has a fixed point $\overline{\chi}$  which contradicts
$q^2=\|\overline{\chi}\|^2=-q\frac{\left<G(\overline{\chi)},\overline{\chi}\right>}{\|G(\overline{\chi})\|} < 0$.

To show the condition needed for $\left<G_h(\chi), \chi\right>$ for $\chi\in V_{h,0}^T$, we  obtain
\begin{align}
\begin{split}
\left<G_h(\chi), \chi\right>
=&
\left<\Lambda_h\chi,\chi\right>-\left<b,\chi\right>
+
\Big<\sin(\frac{\tau\Lambda_h^{\frac12}}2)\widetilde{P}_h\int_0^1f\big(c+2\sin(\frac{\tau\Lambda_h^{\frac12}}2\big) \theta \chi)\,\dd \theta,\chi\Big>
\\
=&
\|\chi\|_1^2-\left<b,\chi\right>
+
\frac12\Big<\int_0^1f\big(c+2\sin(\frac{\tau\Lambda_h^{\frac12}}2) \theta \chi\big)\,\dd \theta,2\sin(\frac{\tau\Lambda_h^{\frac12}}2)\chi\Big>
\\
\geq
&
\lambda_1\|\chi\|^2-\frac{\lambda_1}2\|\chi\|^2-\frac1{2\lambda_1}\|b\|^2
+
\frac12\int_\mathcal{D}F\big(2\sin(\frac{\tau\Lambda_h^{\frac12}}2)\chi+c\big)\,\dd x-\frac12\int_\mathcal{D}F(c)\,\dd x,
\end{split}
\end{align}
which is positive if $\|\chi\| $ is large enough. Therefore,  the problem  \eqref{eq:auxiliary-porblem-I} admits a solution in $B_q$. Let $\chi_1$ and $\chi_2$ be the two
solutions of \eqref{eq:auxiliary-porblem-I} in $V_{h,0}^T$. Then, by one-sided Lipschitz condition \eqref{eq:one-side-L-condition}
\begin{align}
\begin{split}
2\|\Lambda_h^{\frac12}(\chi_1-\chi_2)\|^2
=
&
-\int_0^1\big<f\big(U_h^n+2\sin(\frac{\tau\Lambda_h^{\frac12}}2) \theta \chi_1\big)
-
f\big(U_h^n+2\sin(\frac{\tau\Lambda_h^{\frac12}}2) \theta \chi_2\big),\;2\sin(\frac{\tau\Lambda_h^{\frac12}}2)(\chi_1-\chi_2)\big>\dd \theta
\\
\leq
&
\int_0^1\theta\|2\sin(\frac{\tau\Lambda_h^{\frac12}}2)(\chi_1-\chi_2)\|^2\dd \theta
\\
\leq
&
 \frac{\tau^2}2\|\Lambda_h^{\frac12}(\chi_1-\chi_2)\|^2.
\end{split}
\end{align}
Therefore,  if $\tau<2$,  the problem  \eqref{eq:auxiliary-porblem-I} admits one unique solution
and the proof is complete. $\hfill\square$


In the next lemma, we show the fully discrete problem \eqref{eq:U-formulation}-\eqref{eq:V-formulation} preserves the following energy conservation property.
\begin{lemma}\label{lem:Hamiltonian-full-auxiliary}
The fully discrete problem \eqref{eq:U-formulation}-\eqref{eq:V-formulation}  exactly preserves the following Hamiltonian,
\begin{align}\label{lem:eq-hamiltonian}
J(U_h^{n+1},\overline{V}_h^{n+1})= J(U_h^n,V_h^n), \;n=1,2,\cdots,N,
\end{align}
where $\overline{V}^{n+1}:=V^{n+1}-P_h\delta W^n$.
\end{lemma}
{\it Proof of Lemma \ref{lem:Hamiltonian-full-auxiliary}.} To show this lemma, we recall the fully discrete problem \eqref{eq:U-formulation}-\eqref{eq:V-formulation} and then use the fact $\overline{V}^{n+1}=V^{n+1}-P_h\delta W^n$ to obtain
\begin{align}
U_h^{n+1}=C_h(\tau)U_h^n+\Lambda_h^{-\frac12}S_h(\tau)V_h^n
-\Lambda_h^{-1}(1-C_h(\tau))P_h\nabla_{AVF}f(U_h^n,U_h^{n+1}),\label{eq:U-formulation-I}
\\
\overline{V}_h^{n+1}
=
-\Lambda_h^{\frac12}S_h(\tau)U_h^n+C_h(\tau)V^n_h-\Lambda_h^{-\frac12}S_h(\tau)P_h \nabla_{AVF}f(U_h^n,U_h^{n+1}).
\label{eq:V-formulation-I}
\end{align}
Inserting  \eqref{eq:U-formulation-I} and \eqref{eq:V-formulation-I} into \eqref{eq:definition-J}, we have
\begin{align}\label{eq:J-expresstion}
\begin{split}
J(U_h^{n+1},\overline{V}_h^{n+1})
=
&
\frac12\big<(C^2_h(\tau)+S^2_h(\tau))V_h^n, V_h^n\big>
+
\frac12\big<((\Lambda_h^{\frac12}S_h(\tau))^2+(\Lambda_h^{\frac12}C_h(\tau))^2)U_h^n,U_h^n\big>
\\
&-
\big<\big(C_h(\tau)(1-C_h(\tau))-(S^2_h(\tau))\big)U_h^n, \nabla_{AVF}f(U_h^n,U_h^{n+1})\big>
\\
&
-
\big<\big(\Lambda_h^{-\frac12}S_h(\tau)(I-C_h(\tau))+\Lambda_h^{-\frac12}C_h(\tau)S_h(\tau)\big)V_h^n,
\nabla_{AVF}f(U_h^n,U_h^{n+1})\big>
\\
&
+
\frac12\big<\Lambda_h^{-1}((I-C_h(\tau))^2+S^2_h(\tau))P_h\nabla_{AVF}f(U_h^n,U_h^{n+1}),P_h\nabla_{AVF}f(U_h^n,U_h^{n+1})\big>
\\
&
+
\int_\mathcal{D}F(U_h^{n+1})\,\dd x
\\
=
&
\frac12\|\nabla U_h^n\|^2
+
\frac12\|V_h^n\|^2
+
\big<(I-C_h(\tau))U_h^n,\nabla_{AVF}f(U_h^n,U_h^{n+1})\big>
\\
&
-
\big<\Lambda_h^{-\frac12}S_h(\tau)V_h^n,\nabla_{AVF}f(U_h^n,U_h^{n+1})\big>
\\
&
+
\big<\Lambda_h^{-1}(I-C_h(\tau))P_h\nabla_{AVF}f(U_h^n,U_h^{n+1}),\nabla_{AVF}f(U_h^n,U_h^{n+1})\big>
+
\int_\mathcal{D}F(U_h^{n+1})\,\dd x.
\end{split}
\end{align}
It follows from \eqref{eq:U-formulation} that
\begin{align}
U_h^{n+1}-U_h^n
=
(C_h(\tau)-I)U_h^n+\Lambda_h^{-\frac12}S_h(\tau)V_h^n-\Lambda_h^{-1}(I-C_h(\tau))P_h\nabla_{AVF}f(U_h^n,U_h^{n+1}).
\end{align}
Equation \eqref{eq:J-expresstion} then can be accordingly rewritten as
\begin{align}\label{eq:J-relation-II}
\begin{split}
J(U_h^{n+1},\overline{V}_h^{n+1})
=
&
\frac12\|\nabla U_h^n\|^2
+
\frac12\|V_h^n\|^2
+
\big<U_h^n-U_h^{n+1},P_h\nabla_{AVF}f(U_h^n,U_h^{n+1})\big>
+
\int_\mathcal{D}F(U_h^{n+1})\,\dd x.
\end{split}
\end{align}
By Assumption \ref{assum:Nonlinearity} and the fact $\nabla_{AVF}f\big(U_h^n,U_h^{n+1})=\int_0^1f(U_h^n+\theta(U_h^{n+1}-U_h)\big)\,\dd \theta$, one can infer that
\begin{align}
\int_\mathcal{D}F(U_h^{n+1})\,\dd x-
\int_\mathcal{D}F(U_h^n)\,\dd x
=
\big<\nabla_{AVF}f(U_h^n,U_h^{n+1}), U_h^{n+1}-U_h^n\big>.
\end{align}
Inserting the above formula into \eqref{eq:J-relation-II} implies
\begin{align}
\begin{split}
J(U_h^{n+1},\overline{V}_h^{n+1})
=
&
\frac12\|\nabla U_h^n\|^2
+
\frac12\|V_h^n\|^2
+
\big<U_h^n-U_h^{n+1},P_h\nabla_{AVF}f(U_h^n,U_h^{n+1})\big>
+
\int_\mathcal{D}F(U_h^{n+1})\,\dd x
\\
=
&
\frac12\|\nabla U_h^n\|^2
+
\frac12\|V_h^n\|^+
\int_\mathcal{D}F(U_h^{n})\,\dd x
\\
=
&
J(U_h^n,V_h^n).
\end{split}
\end{align}
This proof is complete. $\hfill\square$

In the following lemma, we will follow the arguments used in the proof of (5.43) in \cite[Section 5.4]{Cox2013Local} and \cite[Proposition 4.2]{cuistrong} to derive the following
exponential integrability property of the numerical of \eqref{eq:U-formulation}-\eqref{eq:V-formulation}.
\begin{lemma}\label{lem:bound-solution-full}
Suppose Assumptions \ref{assum:Nonlinearity}-\ref{assum:initial-value-u0} hold and let
 $X_h^n=\big[U_h^n,V_h^n\big]^\top$ be the solution of the problem \eqref{eq:U-formulation}-\eqref{eq:V-formulation}.
 Then, the following two estimates hold, for any $p\geq 2$
\begin{align}\label{lem:eq-bound-J-full}
\sup_{n\in\{1,2,\cdots,N\}}\mathbb{E}\big[J^p(U_h^n,V_h^n)\big]<\infty,
\end{align}
and
\begin{align}\label{lem:expon-inte-solution-full}
\mathbb{E}\Big(\exp\Big(c\tau\sum_{i=1}^N\|U_h^i\|_{L^6}^2\Big)\Big)<\infty.
\end{align}
\end{lemma}
{\it Proof of Lemma \ref{lem:bound-solution-full}.}
We start by showing \eqref{lem:eq-bound-J-full}.  By \eqref{lem:eq-hamiltonian} and the fact $\overline{V}_h^{j+1}=V_h^{j+1}-P_h\delta W^j$, one can observe that
\begin{align}\label{eq:J(j+1)-relation-J(j)}
\begin{split}
J(U_h^{j+1},V_h^{j+1})
=&
J(U_h^{j+1},\overline{V}_h^{j+1})+\frac12\big<P_h\delta W^j,P_h\delta W^j\big>+\big<\overline{V}_h^{j+1}, P_h\delta W^j \big>
\\
=
&
J(U_h^j,V_h^j)+\frac12\big<P_h\delta W^j,P_h\delta W^j\big>+\big<\overline{V}_h^{j+1}, P_h\delta W^j \big>.
\end{split}
\end{align}
Summing with respect to  $j$ in \eqref{eq:J(j+1)-relation-J(j)}, thus yields
\begin{align}
J(U_h^{n},V_h^{n})
=
J(U_h^0,V_h^0)
+
\sum_{i=0}^{n-1}\left( \frac12\big<P_h\delta W^i,P_h\delta W^i\big>+\big<\overline{V}_h^{i+1}, P_h\delta W^i \big>\right).
\end{align}
We drop the sum on the left, take the $p$th power and the supremum with respect
to $n$, and then the expectation to get
\begin{align}\label{eq:J-bound-full-I}
\begin{split}
\mathbb{E}\sup_{1\leq n\leq N}\big[J^p(U_h^{n},V_h^{n})\big]
\leq
&
C
\mathbb{E}\sup_{1\leq n\leq N}\Big[\Big(J(U_h^0,V_h^0)
+
\sum_{i=0}^{n-1}\|P_h\delta W^i\|^2
+
\sum_{i=0}^{n-1}\big<\overline{V}_h^{i+1}, P_h\delta W^i \big>\Big)^p\Big]
\\
\leq
&
C
\Big(\mathbb{E}[J^p(U_h^0,V_h^0)]
+
\mathbb{E}\Big(\sum_{i=0}^{N-1}\|P_h\delta W^i\|^2\Big)^p
+
\mathbb{E}\sup_{1\leq n\leq N}\Big(\sum_{i=0}^{n-1}\big<\overline{V}_h^{i+1}, P_h\delta W^i \big>\Big)^p\Big).
\end{split}
\end{align}
By the
Burkholder-Davis-Gundy inequality, we have
\begin{align}\label{eq:expection-summ-W}
\begin{split}
\mathbb{E}\Big[\Big(\sum_{i=0}^{N-1}\|P_h\delta W^i\|^2\Big)^p\Big]
\leq
&
CN^{p-1}\sum_{i=0}^{N-1}\mathbb{E}\Big[\|P_h\delta W^i\|^{2p}\Big]
\\
\leq
&
CN^{p-1}\sum_{i=0}^{N-1}\tau^p\Big[\|P_hQ^{\frac12}\|_{\mathcal{L}_2(H)}^{2p}\Big]
\leq
CT^{p-1}\|Q^{\frac12}\|_{_{\mathcal{L}_2(H)}}^{2p}.
\end{split}
\end{align}
Moreover, utilizing Cauchy's inequality and the
Burkholder-Davis-Gundy inequality and applying \eqref{lem:eq-hamiltonian} and the fact that $\overline{V}^{i+1}$ is independent of $\delta W^i$ yields
\begin{align}\label{eq:sum-V-W-expection}
\begin{split}
\mathbb{E}\sup_{1\leq n\leq N}\Big[\Big(\sum_{i=0}^{n-1}(\overline{V}_h^{i+1}, P_h\delta W^i )\Big)^p\Big]
\leq
&
C\mathbb{E}\Big(\sum_{i=0}^{N-1}\tau\|P_hQ^{\frac12}\overline{V}_h^{i+1}\|^2\Big)^{\frac p2}
\\
\leq
&
C\Big(1+\mathbb{E}\Big(\sum_{i=0}^{N-1}\tau\|P_hQ^{\frac12}\overline{V}_h^{i+1}\|^2\Big)^p\Big)
\\
\leq
&
C\Big(1+(N \tau)^{p-1}\|Q^{\frac12}\|^{2p}\sum_{i=0}^{N-1}\tau \mathbb{E} \|\overline{V}_h^{i+1}\|^{2p}\Big)
\\
\leq
&
C\Big(1+(N\tau)^{p-1}\|Q^{\frac12}\|^{2p}\sum_{i=0}^{N-1}\tau \mathbb{E} [J^p(U_h^{i+1},\overline{V}_h^{i+1})]\Big)
\\
=
&
C\Big(1+(N\tau)^{p-1}\|Q^{\frac12}\|^{2p}\sum_{i=0}^{N-1}\tau \mathbb{E} [J^p(U_h^i,V_h^i)]\Big).
\end{split}
\end{align}
As $N\tau=T$ and $\|Q^{\frac12}\|\leq \|Q^{\frac12}\|_{\mathcal{L}_2(H)}$, by inserting \eqref{eq:expection-summ-W} and \eqref{eq:sum-V-W-expection} into
\eqref{eq:J-bound-full-I}
\begin{align}
\begin{split}
\mathbb{E}\sup_{1\leq n\leq N}\big[J^p(U_h^{n},V_h^{n})\big]
\leq
&
C(T,p,\gamma)\Big(1+\sum_{i=0}^{N-1}\tau \mathbb{E} [J^p(U_h^i,V_h^i)]\Big)
\\
\leq
&
C(T,p,\gamma)\Big(1+\sum_{n=0}^{N-1}\tau \mathbb{E}\sup_{0\leq i\leq n} [J^p(U_h^i,V_h^i)]\Big).
\end{split}
\end{align}
Hence  the inequality \eqref{lem:eq-bound-J-full} follows from Gronwall's lemma.

To show \eqref{lem:expon-inte-solution-full}, we introduce the following auxiliary problem,  for $t\in[t_n, t_{n+1}]$
\begin{align}
\begin{split}
\left\{\begin{array}{ll}
\dd U_h(t)=0,& U_h(t_n)=U_h^{n+1},
\\
\dd V_h(t)=P_h \,\dd W(t),& V_h(t_n)=\overline{V}_h^{n+1},
\end{array}\right.
\end{split}
\end{align}
and a linear operator as follows
\begin{align}
\begin{split}
(\mathbb{D}_{\widetilde{\mathbb{F}}_h,\mathbb{B}_h}J)(u_h,v_h):
=
\frac12 \sum_{i=1}^\infty\big<D_{v_h,v_h}J(u_h,v_h)P_hQ^{\frac12}e_i,P_hQ^{\frac12}e_i\big>,
\end{split}
\end{align}
where $\widetilde{\mathbb{F}}_h=\big[0,0\big]^\top$ and $\mathbb{B}_h=\big[0,P_h\big]^\top$.
Adapting similar arguments used in the proof of  (5.43) in \cite[Section 5.4]{Cox2013Local} leads to, for $t\in[t_n,t_{n+1}]$
\begin{align}
\begin{split}
&(\mathbb{D}_{\widetilde{\mathbb{F}}_h,\mathbb{B}_h}J)(U_h(t),V_h(t))
+
\frac1{2\exp(\alpha t)}\sum_{i=1}^\infty\left<(P_hQ^{\frac12})^*V_h(t),e_i\right>^2
\\
=
&
\frac12 \mathrm{Tr}(P_hQP_h)+\frac1{2\exp(\alpha t)}\sum_{i=1}^\infty\left<V_h(t),Q^{\frac12}e_i\right>^2
\\
\leq
&
\frac12 \mathrm{Tr}(Q)+\frac1{\exp(\alpha t)}J(U_h(t),V_h(t))\mathrm{Tr}(Q).
\end{split}
\end{align}
Then, from the exponential integrability lemma in \cite{Cox2013Local}, it follows that, for $t\in[t_n,t_{n+1}]$
\begin{align}
\begin{split}
\mathbb{E}\Big[\exp\Big(\frac{J(U_h(t),V_h(t))}{\exp(\alpha t)}
-
\int_{t_n}^t \frac {\widetilde{U}(s)}{2\exp(\alpha s)}\,\dd s\Big)\Big]
\leq
&
\mathbb{E}\Big[\exp\Big(\frac{J(U_h(t_n),V_h(t_n))}{\exp(\alpha t_n)}\Big)\Big]
\\
=
&\mathbb{E}\Big[\exp\Big(\frac{J(U_h^{n+1},\overline{V}_h^{n+1})}{\exp(\alpha t_n)}\Big)\Big],
\end{split}
\end{align}
where $\alpha>\mathrm{Tr}(Q)$ and $\widetilde{U}=-\frac12\mathrm{Tr}(Q)$.
This together with \eqref{lem:eq-hamiltonian} and the fact
$\big[U_h(t_{n+1}), V_h(t_{n+1})\big]^\top=\big[U_h^{n+1}, V_h^{n+1}\big]^\top$ leads to
\begin{align}
\begin{split}
\mathbb{E}\Big[\exp\Big(\frac{J(U_h^{n+1},V_h^{n+1})}{\exp(\alpha t_{n+1})}\Big)\Big]
\leq
&
\exp\Big(\frac {\tau \mathrm{Tr}(Q)}2\Big)\mathbb{E}\Big[\exp\Big(\frac{J(U_h^{n},V_h^{n})}{\exp(\alpha t_n)}\Big)\Big]
\\
\leq
&
\exp\Big(\frac {t_{n+1} \mathrm{Tr}(Q)}2\Big)\mathbb{E}\Big[\exp\Big(J(U_h^0,V_h^0)\Big)\Big]
<
\infty.
\end{split}
\end{align}
Therefore, this in combination with  Jensen's inequality, Young's inequality and the fact $\|u_h\|_{L^6}\leq c\|u_h\|^{a}\|\Lambda_h^{\frac12}u_h\|^{1-a}$ with $a=\frac d3$, $d=1,2$  implies
\begin{align}
\begin{split}
\mathbb{E}\Big[\exp\Big(C\tau\sum_{i=0}^n\|U_h^i\|_{L^6}^2\,\dd s\Big)\Big]
\leq
&
\sup_{i\in\{1,2,\cdots,N\}}\mathbb{E}\big[\exp(CT\|U_h^i\|_{L^6}^2)\big]
\\
\leq
&
\sup_{i\in\{1,2,\cdots,N\}}\mathbb{E}\Big[\exp\Big(\frac{\|\nabla U_h^i\|^2}{2\exp(\alpha t_i)}\Big) \exp\Big(\exp(\frac{a}{1-a}\alpha T)\|U_h^i\|^2(CT)^{\frac1{1-a}}2^{\frac a{1-a}}\Big)\Big]
\\
\leq
&
\sup_{i\in\{1,2,\cdots,N\}}\mathbb{E}\Big[\exp\Big(\frac{\| \nabla U_h^i\|^2}{2\exp(\alpha t_i)}\Big) \exp\Big(\frac{\|U_h^i\|_{L^4}^4}{2\exp(\alpha t_i)}+C(T)\Big)\Big]
\\
\leq
&
C
(T)\sup_{i\in\{1,2,\cdots,N\}}\mathbb{E}\Big[\exp\Big(\frac{J(U_h^i,V_h^i)}{2\exp(\alpha t_i)}\Big)\Big]
<
\infty.
\end{split}
\end{align}
Therefore, the desired result is shown and the proof is complete. $\hfill\square$

Similar to the semi-discrete case, the exact solution of the fully discrete problem \eqref{eq:solution-full} also satisfies the following trace formula.
\begin{lemma}\label{eq-lem:discrete-trace-formula-full}
Let Assumptions \ref{assum:Nonlinearity}-\ref{assum:initial-value-u0} be fulfilled and $X_h^n=\big[U_h^n,V_h^n\big]^\top$ be the solution of the problems \eqref{eq:U-formulation}-\eqref{eq:V-formulation}. Then $X_h^n$ satisfies the trace formula
\begin{align}\label{lem:discrete-trace-formula-full}
\mathbb{E}[J(U_h^n,V_h^n)]=\mathbb{E}[J(U_h^0,V_h^0)]+\frac12 \mathrm{Tr}(P_hQP_h)t_n.
\end{align}
\end{lemma}

\subsection{Strong convergence rates of the fully discretization}

The forthcoming  lemma, quoted from \cite{anton2016full}, is a spatial version of Lemma
\ref{lem:eq-temporal-result-E}, and plays a significant role in the error estimates of the fully discrete approximation.
\begin{lemma}
For any $\alpha\in[0,1]$, there exists a constant $C=C(r)$ such that
\begin{align}\label{eq:temporal-regulari-s-c-sem}
\|(S_h(t)-S_h(s))\Lambda_h^{-\frac \alpha2}\|_{\mathcal{L}(H)}
\leq
C (t-s)^\alpha,\quad  \|(C_h(t)-C_h(s))\Lambda_h^{-\frac \alpha2}P_h\|_{\mathcal{L}(H)}
\leq
C (t-s)^\alpha,
\end{align}
and
\begin{align}\label{eq:temporal-Eh}
\|(E_h(t)-E_h(s))\chi_h\|_{\mathbb{H}^0}
\leq
C(t-s)^\alpha\|\chi_h\|_{\mathbb{H}^\alpha},\; \forall \chi_h=[\chi_h^1,\chi_h^2]^\top\in V_h\times V_h
\end{align}
for all $t\geq s\geq 0$.
\end{lemma}

\begin{theorem}\label{them:covergence-full}
Let  Assumptions \ref{assum:Nonlinearity}-\ref{assum:initial-value-u0}, \ref{eq:assumption-rh-ph} be fulfilled and
$X_h^n=[U_h^n,V_h^n]^\top$ be the solution of the scheme \eqref{eq:U-formulation}-\eqref{eq:V-formulation}.
Then, for all $n\in\{1,2,\cdots, N\}$, it holds
\begin{align}\label{them-eq:error-full-u}
\|X(t_n)-X_h^n\|_{L^p(\Omega;\mathbb{H}^0)}
\leq
C(h^{\frac{r+1}{r+2}\min\{\gamma,r+2\}}+\tau).
\end{align}
Additionally, for the velocity  in  dimension one,  it holds for all $n\in\{1,2,\cdots, N\}, N\in \mathbb{N}^+$,
\begin{align}\label{them-eq:error-full-v}
\|v(t_n)-V_h^n\|_{L^p(\Omega;H)}
\leq
C(h^{\frac{r+1}{r+2}\min\{\gamma-1,r+2\}}+\tau^{\min\{\gamma-1,1\}}).
\end{align}
Furthermore, if $\Lambda_hP=P_h\Lambda_h$, it holds, for all $n\in\{1,2,\cdots, N\}, N\in \mathbb{N}^+$,
\begin{align}\label{them-eq:error-full-u-ii}
\|X(t_n)-X_h^n\|_{L^p(\Omega;\mathbb{H}^0)}
\leq
C(h^{\min\{\gamma,r+1\}}+\tau),
\end{align}
and for the velocity  in  dimension one,
\begin{align}\label{them-eq:error-full-v-ii}
\|v(t_n)-V_h^n\|_{L^p(\Omega;H)}
\leq
C(h^{\min\{\gamma-1,r+1\}}+\tau^{\min\{\gamma-1,1\}}).
\end{align}

\begin{remark}\label{eq:error-full-Galerkin-spetral}
For the spectral Galerkin method, Assumption \ref{eq:assumption-rh-ph}  holds for every $r\in \mathbb{N}^+$ and $\Lambda_hP_h=P_h\Lambda_h$. Therefore, the errors of the spectral Galerkin method can be estimated as  follows,
\begin{align}\label{them-eq:error-semi-norm-H-ii-i}
\|X(t_n)-X_h^n\|_{L^p(\Omega;\mathbb{H}^0)}
\leq
C ( h^{\gamma}+\tau),\;with\;h=\lambda_{J+1}^{-\frac12},
\end{align}
and for the velocity in  dimension one
\begin{align}\label{them-eq:vecolcity-error-semi-norm-H-ii-ii}
\|v(t)-v_h(t)\|_{L^p(\Omega;H)}
\leq
C(\lambda_{J+1}^{\gamma-1}+\tau^{\min\{\gamma-1,1\}}),\;with\;h=\lambda_{J+1}^{-\frac12}.
\end{align}
\end{remark}

\end{theorem}
{\it Proof of Theorem \ref{them:covergence-full}.}
Equivalently, \eqref{eq:solution-full} can be reformulated as
\begin{align}
X_h^n
=
E_h(t_n)P_hX_0
+
\sum_{i=0}^{n-1}
\int_{t_i}^{t_{i+1}} E_h(t_n-s)\mathbb{F}_h^i\,\dd s
+
\sum_{i=0}^{n-1}\int_{t_i}^{t_{i+1}} E_h(t_{n-i-1})\mathbb{B}_h\,\dd W(s).
\end{align}
Therefore, the difference between $X_h^n$ and $X(t_n)$ can be decomposed as follows:
\begin{align}
\begin{split}
X(t_n)-X^n_h
=
&
(E(t_n)-E_h(t_n)P_h)X_0
\\
&
+
\sum_{i=0}^{n-1}\int_{t_i}^{t_{i+1}}\big(E(t_n-s)\mathbb{F}(X(s))
-
E_h(t_n-s)\mathbb{F}_h^i\big)\,\,\dd s
\\
&
+
\sum_{i=0}^{n-1}\int_{t_i}^{t_{i+1}}\big(E(t_n-s)\mathbb{B}-E_h(t_{n-1-i})\mathbb{B}_h\big) \,\dd W(s)\\
=&\mathbb{I}+\mathbb{II}+\mathbb{III}.
\end{split}
\end{align}
In the same manner as \eqref{eq:estimation-I1}, the first term $\mathbb{I}_1$ can be estimated as follows
\begin{align}
\|\mathbb{I}\|_{\mathbb{H}^0}
\leq
Ch^{\frac{r+1}{r+2}\min\{\gamma,r+2\}}\|X_0\|_{\mathbb{H}^{\min\{\gamma,r+1\}}}.
\end{align}
In order to deal with the term $\mathbb{II}$, we decompose it into four terms as follows:
\begin{align}
\begin{split}
\mathbb{II}
=
&
\sum_{i=0}^{n-1}\int_{t_i}^{t_{i+1}}\big(E(t_n-s)-E_h(t_n-s)P_h\big)\mathbb{F}(X(s))\,\dd s
\\
&+
\sum_{i=0}^{n-1}\int_{t_i}^{t_{i+1}}E_h(t_n-s)\big(\mathbb{F}_h(X(s))-\mathbb{F}_h(X(t_i))\big)\,\dd s
\\
&+
\sum_{i=0}^{n-1}\int_{t_i}^{t_{i+1}}E_h(t_n-s)\big(\mathbb{F}_h(X(t_i))-\mathbb{F}_h(X_h^i)\big)\,\dd s
\\
&+
\sum_{i=0}^{n-1}\int_{t_i}^{t_{i+1}}E_h(t_n-s)\big(\mathbb{F}_h(X_h^i)-\mathbb{F}_h^i\big)\,\dd s
\\
=&
\mathbb{II}_1+\mathbb{II}_2+\mathbb{II}_3+\mathbb{II}_4.
\end{split}
\end{align}
In the sequel we treat with the above four terms one by one. Similar  to \eqref{eq:estimation-I2}, we derive
\begin{align}
\begin{split}
\|\mathbb{II}_1\|_{\mathbb{H}^0}
\leq
Ch^{\frac{r+1}{r+2}\min\{\gamma,r+2\}}\sum_{i=0}^{n-1}\int_{t_i}^{t_{i+1}}\|\Lambda^{\frac{\gamma-1}2}f(u(s))\|\,\dd s.
\end{split}
\end{align}
For the term $\mathbb{II}_2$, using  \eqref{eq:relation-H(-1)-L78} and the boundness of
 the cosine and sine operators implies
\begin{align}\label{eq:bound-II2}
\begin{split}
\|\mathbb{II}_2\|_{\mathbb{H}^0}
\leq&
\sum_{i=0}^{n-1}\int_{t_i}^{t_{i+1}}\|E_h(t_n-s)(\mathbb{F}_h(X(s))-\mathbb{F}_h(X(t_i)))\|_{\mathbb{H}^0}\,\dd s
\\
\leq
&
C\sum_{i=0}^{n-1}\int_{t_i}^{t_{i+1}}\|\Lambda_h^{-\frac12}P_h(f(u(s))-f(u(t_i)))\|\,\dd s
\\
\leq
&
C\sum_{i=0}^{n-1}\int_{t_i}^{t_{i+1}}\|f(u(s))-f(u(t_i))\|_{L^{\frac65}}\,\dd s
\\
\leq
&
C\sum_{i=0}^{n-1}\int_{t_i}^{t_{i+1}}\|u(s)-u(t_i)\|(1+\|u(s)\|_{L^6}^2+\|u(t_i)\|_{L^6}^2)\,\dd s.
\end{split}
\end{align}
 Similarly as above, we obtain
\begin{align}
\begin{split}
\|\mathbb{II}_3\|_{\mathbb{H}^0}
\leq
&
\sum_{i=0}^{n-1}\int_{t_i}^{t_{i+1}}\|E_h(t_n-s)P_h(\mathbb{F}(X(t_i))-\mathbb{F}(X_h^i))\|_{\mathbb{H}^0}\,\dd s
\\
\leq
&
C\sum_{i=0}^{n-1}\tau\|u(t_i)-U_h^i\|(1+\|u(t_i)\|^2_{L^6}+\|U_h^i\|^2_{L^6}).
\end{split}
\end{align}
To bound the term $\mathbb{II}_4$, we first employ \eqref{eq:U-formulation}, \eqref{eq:temporal-regulari-s-c-sem} and the fact $C_h(0)=P_h$ and $S_h(0)=0$ to derive
\begin{align}\label{eq:error-u-n-u-n+1}
\begin{split}
\|U_h^{n+1}-U_h^n\|
\leq
&
\|(C_h(\tau)-C_h(0))U_h^n\|
+
\|\Lambda_h^{-\frac12}(S_h(\tau)-S_h(0))V_h^n\|
\\
&+
\Big\|\Lambda_h^{-1}(C_h(0)-C_h(\tau))P_h\int_0^1f(U^n_h+\theta(U^{n+1}_h-U^n_h))\,\dd \theta\Big\|
\\
\leq
&
C\tau\|\Lambda_h^{\frac12}U_h^n\|
+
C\tau\|V_h^n\|
\\
&+
C\tau\Big\|\Lambda_h^{-\frac12}P_h\int_0^1f\big(U^n_h+\theta(U^{n+1}_h-U^n_h)\big)\,\dd \theta\Big\|
\\
\leq
&
C\tau\|\Lambda_h^{\frac12}U_h^n\|
+
C\tau\|V_h^n\|
\\
&+
C\tau\int_0^1\|f\big(U^n_h+\theta(U^{n+1}_h-U^n_h)\big)\|_{L^{\frac65}}\,\dd \theta
\\
\leq
&
C\tau\|\Lambda_h^{\frac12}U_h^n\|
+
C\tau\|V_h^n\|
+
C\tau(1+\|U_h^{n+1}\|_{L^6}^3+\|U_h^n\|^3_{L^6})
\\
\leq
&
C\tau(1+\|\nabla U_h^n\|^3
+
\|V_h^n\|
+
\|\nabla U_h^{n+1}\|^3).
\end{split}
\end{align}
Then the above estimate together with  the similar arguments used in  \eqref{eq:bound-II2} implies
\begin{align}
\begin{split}
\|\mathbb{II}_4\|_{\mathbb{H}^0}
\leq
&
\sum_{i=0}^{n-1}\int_{t_i}^{t_{i+1}}\|E_h(t_n-s)P_h(\mathbb{F}(X_h^i)-\mathbb{F}^i)\|_{\mathbb{H}^0}\,\dd s
\\
\leq&
\sum_{i=0}^{n-1}\int_{t_i}^{t_{i+1}}\left\|\Lambda_h^{-\frac12}P_h\int_0^1f(U^i_h)-f(U_h^i+\theta(U_h^{i+1}-U_h^i))\dd \theta\right\|\,\dd s
\\
\leq
&
C\sum_{i=0}^{n-1}\tau\|U_h^{i+1}-U_h^i\|(1+\|U_h^i\|^2_{L^6}+\|U_h^{i+1}\|^2_{L^6})
\\
\leq
&
C\tau^2\sum_{i=0}^{n-1}\big(1+\|\nabla U_h^i\|^3
+
\|V_h^i\|
+
\|\nabla U_h^{i+1}\|^3\big)\big(1+\|U_h^i\|^2_{L^6}+\|U_h^{i+1}\|^2_{L^6}\big)
\\
\leq
&
C\tau^2\sum_{i=0}^{n-1}(1+\|U_h^i\|_1^5+\|U_h^{i+1}\|_1^5+\|V_h^i\|^2).
\end{split}
\end{align}
Therefore, gathering the above estimates together gives
\begin{align}
\|X(t_n)-X_h^n\|_{\mathbb{H}^0}
\leq
C\sum_{i=0}^{n-1}\|X(t_i)-X_h^i\|_{\mathbb{H}^0}\Psi_i
+
C\sum_{i=0}^{n-1}\Phi_i
+
\|\mathbb{III}\|_{\mathbb{H}^0}
+
Ch^{\frac{r+1}{r+2}\min\{\gamma,r+2\}}\|X_0\|_{\mathbb{H}^\gamma},
\end{align}
where $\Psi_i=C\tau(1+\|U_h^i\|_{L^6}^2+\|u(t_i)\|_{L^6}^2)$ and
\begin{align}
\begin{split}
\Phi_i
=
&
Ch^{\frac{r+1}{r+2}\min\{\gamma,r+2\}}\int_{t_i}^{t_{i+1}}\|\Lambda^{\frac{\gamma-1}2}f(u(s))\|\,\dd s
\\
&
+
C\int_{t_i}^{t_{i+1}}\|u(s)-u(t_i)\|(1+\|u(s)\|_{L^6}^2+\|u(t_i)\|_{L^6}^2)\,\dd s
\\
&
+
C\tau^2(1+\|U_h^i\|_1^5+\|U_h^{i+1}\|_1^5+\|V_h^i\|^2).
\end{split}
\end{align}
By Gronwall's inequality, we have
\begin{align}
\|X(t_n)-X_h^n\|_{\mathbb{H}^0}
\leq
C\Big(\sum_{i=0}^{n-1}\Phi_i
+
\|\mathbb{III}\|_{\mathbb{H}^0}
+
h^{\frac{r+1}{r+2}\min\{\gamma,r+2\}}\|X_0\|_{\mathbb{H}^\gamma}\Big)
\exp\Big(\sum_{i=0}^{n-1}\Psi_i
\Big).
\end{align}
Taking the $p'$th power and  the expectation, and then using H\"{o}lder inequality and the exponential integrability properties of $u(t_n)$ and $U_h^n$ give
\begin{align}
\begin{split}
\mathbb{E}[\|X(t_n)-X_h^n\|_{\mathbb{H}^0}^p]
\leq
&
C\Big(\mathbb{E}\Big[
\sum_{i=0}^{n-1}\Phi_i
+
\|\mathbb{III}\|_{\mathbb{H}^0}
+
h^{\frac{r+1}{r+2}\min\{\gamma,r+2\}}\|X_0\|_{\mathbb{H}^\gamma}\Big]^{2p}\Big)^{\frac12}
\Big(\mathbb{E}\Big(\exp(2p\sum_{i=0}^{n-1}\Psi_i
\Big)\Big)^{\frac12}
\\
\leq
&
C\Big(\mathbb{E}\Big[
\sum_{i=0}^{n-1}\Phi_i\Big]^{2p}\Big)^{\frac12}
+
C\big(\mathbb{E}\|\mathbb{III}\|^{2p}_{\mathbb{H}^0}\big)^{\frac12}
+
Ch^{\frac{r+1}{r+2}p\min\{\gamma,r+2\}}\big(\mathbb{E}\big[\|X_0\|^{2p}_{\mathbb{H}^\gamma}\big]\big)^{\frac12}
\\
\leq
&
Cn^{p-\frac12}\Big(
\sum_{i=0}^{n-1}\mathbb{E}[\Phi_i]^{2p}\Big)^{\frac12}
+
C\big(\mathbb{E}\|\mathbb{III}\|^{2p}_{\mathbb{H}^0}\big)^{\frac12}
+
Ch^{\frac{r+1}{r+2}p\min\{\gamma,r+2\}}\big(\mathbb{E}\big[\|X_0\|^{2p}_{\mathbb{H}^\gamma}\big]\big)^{\frac12}.
\end{split}
\end{align}
From \eqref{eq:error-semigroup}, \eqref{eq:temporal-Eh} and Burkholder-Davis-Gundy inequality, it follows that
\begin{align}\label{eq:estimate-III}
\begin{split}
\mathbb{E}[\|\mathbb{III}\|^{2p}_{\mathbb{H}^0}]
\leq
&
\mathbb{E}\Big[\Big\|\sum_{i=0}^{n-1}\int_{t_i}^{t_{i+1}}(E(t_n-s)-E_h(t_{n-1-i})P_h)\mathbb{B}\,\dd W(s)\Big\|^{2p}_{\mathbb{H}^0}\Big]
\\
\leq
&
C\Big[\sum_{i=0}^{n-1}\int_{t_i}^{t_{i+1}}\Big\|(E(t_n-s)
-
E_h(t_n-s)P_h)\mathbb{B}Q^{\frac12}\Big\|_{\mathcal{L}_2(\mathbb{H}^0,\mathbb{H}^0)}^2\,\dd s\Big]^p
\\
&
+
C\Big[\sum_{i=0}^{n-1}\int_{t_i}^{t_{i+1}}\Big\|(E_h(t_n-s)
-
E_h(t_{n-1-i}))P_h\mathbb{B}Q^{\frac12}\Big\|_{\mathcal{L}_2(\mathbb{H}^0,\mathbb{H}^0)}^2\,\dd s\Big]^p
\\
\leq
&
Ch^{\frac{r+1}{r+2}2p\min\{\gamma,r+2\}}\|\Lambda^{\frac{\gamma-1}2} Q^{\frac12}\|_{\mathcal{L}_2(H)}^{2p}+
C\tau^{2p}\|Q^{\frac12}\|^{2p}_{\mathcal{L}_2(H)}.
\end{split}
\end{align}
Thanks to \eqref{eq-lem:spatial-regularity-f(u)}, \eqref{lem:eq-bound-J-full}, \eqref{eq-them:spatial-regularity-mild} and \eqref{them-eq:temporal-regularity-u}, we obtain
\begin{align}
\begin{split}
\sum_{j=0}^{n-1}\mathbb{E}[\Phi_i]^{2p}
\leq
&
Ch^{\frac{r+1}{r+2}2p\min\{\gamma,r+2\}}
\sum_{i=0}^{n-1}\mathbb{E}\Big[\int_{t_i}^{t_{i+1}}\|\Lambda^{\frac{\gamma-1}2}f(u(s))\|\,\dd s\Big]^{2p}
\\
&+
C\sum_{i=0}^{n-1}\mathbb{E}\Big[\int_{t_i}^{t_{i+1}}\|u(s)-u(t_i)\|(1+\|u(s)\|_{L^6}^2
+\|u(t_i)\|_{L^6}^2)\,\dd s\Big]^{2p}
\\
&
+
C\sum_{i=0}^{n-1}\mathbb{E}[\tau^{4p}(1+\|U_h^i\|_1^5+\|U_h^{i+1}\|_1^5+\|V_h^i\|^2)^{2p}]
\\
\leq
&
Ch^{\frac{r+1}{r+2}2p\min\{\gamma,r+2\}}\tau^{2p-1}\int_0^{t_n}\sup_{s\in[0,T]}
\mathbb{E}\big[\|\Lambda^{\frac{\gamma-1}2}f(u(s))\|^{2p}\big]\,\dd s
\\
&
+
C\tau^{2p-1}\sum_{i=0}^{n-1}\int_{t_i}^{t_{i+1}}\mathbb{E}
\Big[\|u(s)-u(t_i)\|^{2p}(1+\|u(s)\|^2_{L^6}+\|u(t_i)\|_{L^6}^2)^{2p}\Big]\,\dd s
\\
&
+
C\tau^{4p}\sum_{i=0}^{n-1}\sup_{i\in\{1,2,\cdots,N\}}\mathbb{E}\big[(1+\|U_h^i\|_1^5+\|U_h^{i+1}\|_1^5+\|V_h^i\|^2\big)^{2p}\big]
\\
\leq
&
C\tau^{2p-1}\sum_{i=0}^{n-1}\int_{t_i}^{t_{i+1}}\big(\mathbb{E}[\|u(s)-u(t_i)\|^{4p}\big]\mathbb{E}\big[(1+\|u(s)\|_{L^6}^2+\|u(t_i)\|_{L^6}^2)^{4p}\big]\big)^{\frac12}\,\dd s
\\
&+
Ch^{\frac{r+1}{r+2}2p\min\{\gamma,r+2\}}\tau^{2p-1}
+
C\tau^{4p-1}
\\
\leq
&
Ch^{\frac{r+1}{r+2}2p\min\{\gamma,r+2\}}\tau^{2p-1}
+
C\tau^{4p-1}.
\end{split}
\end{align}
This together with \eqref{eq:estimate-III} yields
\begin{align}
\begin{split}
\mathbb{E}[\|X(t_n)-X_h^n\|^p_{\mathbb{H}^0}]
\leq
&
Ch^{\frac{r+1}{r+2}p\min\{\gamma,r+2\}}+Ch^{\frac{r+1}{r+2}p\min\{\gamma,r+2\}}(n\tau)^{p-\frac12}+C\tau^p(n\tau)^{p-\frac12}
\\
\leq
&
 C(h^{\frac{r+1}{r+2}p\min\{\gamma,r+2\}}+\tau^p).
 \end{split}
\end{align}

We are now in the position to verify \eqref{them-eq:error-full-v}. Similarly as before, the error $\|v(t_n)-V_h^n\|_{L^p(\Omega;H)}$ can be decomposed as follows
\begin{align}
\begin{split}
\|v(t_n)-v^n_h\|_{L^p(\Omega;H)}
\leq
&
\|P_2\big(E(t_n)-E_h(t_n)P_h\big)X_0\|_{L^p(\Omega;H)}
\\
&
+
\sum_{i=0}^{n-1}\int_{t_i}^{t_{i+1}}\|P_2\big(E(t_n-s)\mathbb{F}(X(s))
-
E_h(t_n-s)\mathbb{F}_h^i\big)\|_{L^p(\Omega;H)}\,\,\dd s
\\
&
+
\Big\|\sum_{i=0}^{n-1}\int_{t_i}^{t_{i+1}}P_2\big(E(t_n-s)-E_h(t_{n-1-i})P_h\big)\mathbb{B} \,\dd W(s)\Big\|_{L^p(\Omega;H)}\\
=&\mathbb{L}_1+\mathbb{L}_2+\mathbb{L}_3.
\end{split}
\end{align}
By  \eqref{eq:estimate-L1}, the term $\mathbb{L}_1$ can be estimated as follows
\begin{align}\label{eq:estiamte-L1-full}
\mathbb{L}_1
\leq
Ch^{\frac{r+1}{r+2}\min\{\gamma-1,r+2\}}\|X_0\|_{L^p(\Omega;\mathbb{H}^\gamma)}.
\end{align}
In order to properly handle $\mathbb{L}_2$, we need its further decomposition as follows
\begin{align}
\begin{split}
\mathbb{L}_2
\leq&
\int_0^{t_n}\|P_2\big(E(t_n-s)-E_h(t-s)P_h\big)\mathbb{F}\big(X(s)\big)\|_{L^p(\Omega;H)}\,\,\dd s
\\
&+
\sum_{i=0}^{n-1}\int_{t_i}^{t_{i+1}}\| P_2E_h(t_n-s)\big(\mathbb{F}_h(X(s))
-\mathbb{F}_h(X(t_i))\big)\|_{L^p(\Omega;H)}\,\,\dd s
\\
&+
\sum_{i=0}^{n-1}\int_{t_i}^{t_{i+1}}\| P_2E_h(t_n-s)\big(
\mathbb{F}_h(X(t_i))-\mathbb{F}_h^i\big)\|_{L^p(\Omega;H)}\,\,\dd s
\\
=
&
\mathbb{L}_{21}+\mathbb{L}_{22}+\mathbb{L}_{23}.
\end{split}
\end{align}
By \eqref{eq:estimate-L21}, we have
\begin{align}
\begin{split}
\mathbb{L}_{21}\leq Ch^{\frac{r+1}{r+2}
\min\{\gamma-1,r+2\}}.
\end{split}
\end{align}
Similar to \eqref{eq:estimate-L22}, using \eqref{them-eq:temporal-regularity-u}, \eqref{eq:embedding-equatlity-I}, \eqref{eq:bound-u-H1} and the fact $ \dot{H}^1\subset C(\mathcal{D};\mathbb{R})$ in  dimension one yields
\begin{align}
\begin{split}
\mathbb{L}_{22}
\leq
&
C\sum_{i=0}^{n-1}\int_{t_i}^{t_{i+1}} \|f(u(s))-f(u(t_i))\|_{L^p(\Omega;H)}
\,\dd s
\\
\leq
&
C\sum_{i=0}^{n-1}\int_{t_i}^{t_{i+1}} \|u(s)-u(t_i)\|_{L^{2p}(\Omega;H)}(1+\|u(s)\|^2_{L^{4p}(\Omega;C(\mathcal{D};\mathbb{R}))}+\|u(t_i)\|^2_{L^{4p}(\Omega;C(\mathcal{D};\mathbb{R}))})
\,\dd s
\\
\leq
&
C\tau T (1+\sup_{s\in[0,T]}\|u(s)\|^2_{L^{4p}(\Omega;\dot{H}^1)})
\\
\leq&
 C\tau.
 \end{split}
\end{align}
Similarly as before, we utilize \eqref{them-eq:error-full-u}, \eqref{eq:error-u-n-u-n+1}, \eqref{eq:embedding-equatlity-I}, \eqref{eq:bound-u-H1} and \eqref{lem:eq-bound-J-full},  and apply the boundness of
 the cosine operator and  the fact $\mathbb{F}^i=[0,\int_0^1f(U_h^i+\theta(U_h^{i+1}-U_h^i))\,\dd \theta]'$ to bound $\mathbb{L}_{23}$ as follows
\begin{align}
\begin{split}
\mathbb{L}_{23}
\leq
&
\sum_{i=0}^{n-1}\int_{t_i}^{t_{i+1}} \int_0^1\|f(u(t_i))
-
f\big(U_h^i+\theta(U_h^{i+1}-U_h^i)\big)\|_{L^p(\Omega;H)}\,\dd \theta
\,\,\dd s
\\
\leq
&
\sum_{i=0}^{n-1}\int_{t_i}^{t_{i+1}} (\|u(t_i)-U_h^i\|_{L^{2p}(\Omega;H)}
+
\|U_h^{i+1}-U_h^i\|_{L^{2p}(\Omega;H)})
\\
&
\cdot
(1+\|u(t_i)\|_{L^{4p}(\Omega;C(\mathcal{D};\mathbb{R}))}^2+\|U_h^{i+1}\|_{L^{4p}(\Omega;C(\mathcal{D};\mathbb{R}))}^2
+
\|U_h^i\|_{L^{4p}(\Omega;C(\mathcal{D};\mathbb{R}))}^2)
\\
\leq
&
C(h^{\frac{r+1}{r+2}\min\{\gamma,r+2\}}+\tau)
(1+\sup_{s\in[0,T]}\|u(s)\|_{L^{4p}(\Omega;\dot{H}^1)}^2
+
\sup_{i\in\{1,2,\cdots,N\}}\|U_h^i\|_{L^{4p}(\Omega;\dot{H}^1)}^2)
\\
\leq
&
C(h^{\frac{r+1}{r+2}\min\{\gamma,r+2\}}+\tau)
.
\end{split}
\end{align}
Putting the above three estimates together ensures
\begin{align}\label{eq:L2-esitmate-full}
\mathbb{L}_2
\leq
C(h^{\frac{r+1}{r+2}\min\{\gamma-1,r+2\}}+\tau).
\end{align}
For the term $\mathbb{L}_3$, by \eqref{eq:determinisitic-error-velocity} with $\beta=\min\{\gamma,r+3\}$,
\eqref{eq:temporal-regulari-s-c-sem} with $\alpha=\min\{\gamma-1,1\}$ and  the Burkholder-Davis-Gundy inequality, we obtain
\begin{align}\label{eq:L3-esitmate-full}
\begin{split}
\mathbb{L}_3
\leq
&
\Big(\sum_{i=0}^{n-1}\int_{t_i}^{t_{i+1}}\|P_2\big(E(t_n-s)-E_h(t_{n-1-i})P_h\big)\mathbb{B}Q^{\frac12}\|^2_{\mathcal{L}_2(H)} \,\dd s\Big)^{\frac12}
\\
\leq
&
\Big(\sum_{i=0}^{n-1}\int_{t_i}^{t_{i+1}}\|\big(C(t_n-s)-C_h(t_n-s)P_h\big)Q^{\frac12}\|^2_{\mathcal{L}_2(H)}  \,\dd s\Big)^{\frac12}
\\
&
+
\Big(\sum_{i=0}^{n-1}\int_{t_i}^{t_{i+1}}\|\big(C_h(t_n-s)-C_h(t_{n-1-i})\big)P_hQ^{\frac12}\|^2_{\mathcal{L}_2(H)}  \,\dd s\Big)^{\frac12}
\\
\leq
&
Ch^{\frac{r+1}{r+2}\min\{\gamma-1,r+2\}}\|\Lambda^{\frac{\min\{\gamma,r+3\}-1}2}Q^{\frac12}\|_{\mathcal{L}_2(H)}
+
C\tau^{\min\{\gamma-1,1\}}\|\Lambda^{\frac{\min\{\gamma-1,1\}}2}Q^{\frac12}\|_{\mathcal{L}_2(H)}
\\
\leq
&C(h^{\frac{r+1}{r+2}\min\{\gamma-1,r+2\}}
+
\tau^{\min\{\gamma-1,1\}}).
\end{split}
\end{align}
Then summing up \eqref{eq:estiamte-L1-full}, \eqref{eq:L2-esitmate-full} and \eqref{eq:L3-esitmate-full} shows \eqref{them-eq:error-full-v}. Hence the proof is complete.
$\hfill\square$

\section{Numerical experiments}
In this section, we perform numerical experiments to verify the previous theoretical finding. For simplicity, we consider the stochastic wave equation in one dimension as follows:
\begin{align}\label{eq:sswe-one-dimension}
\begin{split}
\left\{\begin{array}{ll}
\dd u(t) =v(t)\dd t,&(t,x)\in [0, 1]\times [0,1],
\\
\dd v(t)=  u_{xx}(t) \dd t
-u^3 \dd t + \dd W(t),& (t,x)\in [0, 1]\times [0,1],
\\
u(0)=0,\;v(0)=0,&x\in \mathcal{D},
\end{array}\right.
\end{split}
\end{align}
where $\{W(t)\}_{t\in[0,1]}$
stands for a $Q$-Wiener process with the covariance
operator $Q=\Lambda^{-0.5005}$. One can easily see that Assumption 2.3 is fulfilled with $\gamma =1$ for $ Q=\Lambda^{-0.5005}$. In the following experiments, we aim to confirm  the energy-preserving property
of the proposed fully-discrete schemes \eqref{eq:U-formulation}-\eqref{eq:V-formulation} and test mean-square
approximations of the exact solution $(u(t), v(t))$ to \eqref{eq:sswe-one-dimension} at the endpoint $T = 1$. To do this, we take a piece-wise linear continuous finite element method for the spatial discretization and
the proposed time-stepping scheme for the temporal discretization. The true solutions $(u(t),v(t))$ are not available and are computed by numerical ones with small step-sizes $h_{exact}$ and $\tau_{exact}$. The expectation is approximated by computing averages over 1000 samples.

Now let us begin with tests on the   convergence rates in time. The true solutions $(u,v)$ are computed by using $h_{exact}=2^{-8}$ and $\tau_{exact}=2^{-8}$. Meanwhile, numerical simulations are performed with five different temporal step sizes $\tau=2^{-i}, i=2,3,4,5,6$ and  the resulting mean-square errors are depicted in the left hand side of  Figure 1. One can observe the expected convergence rate of order $1$, which coincides with our previous theoretical findings in Theorem \ref{them:covergence-full} and Remark \ref{eq:error-full-Galerkin-spetral}.


To test the convergence rates in space, we also do numerical approximation with five different space-step size $h=2^{-i}$, $i=2,3,4,5,6$. The true solution are approximated by the method \eqref{eq:U-formulation}-\eqref{eq:V-formulation}
with $k_{exact}=2^{-6}$ and $h_{exact}=2^{-9}$. In the right hand side of Figure 1, we present the resulting errors of the proposed method \eqref{eq:U-formulation}-\eqref{eq:V-formulation} in spatial direction. As expected, the numerical performance is all consistent with the previous theoretical results.
\begin{figure}[!ht]\label{fig:space-results}
\centering
 \scalebox{0.5}[0.5]{\includegraphics[0.3in,
0.2in][5.3in,4.1in]{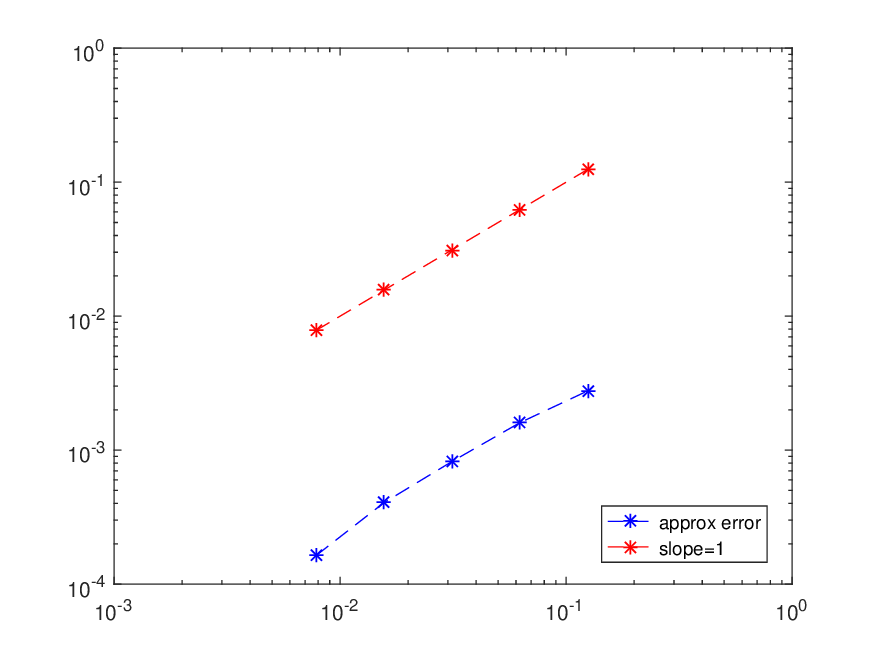}}\qquad\qquad\qquad\qquad\qquad\qquad
 \scalebox{0.5}[0.5]{\includegraphics[0.4in,
0.1in][5.3in,4.1in]{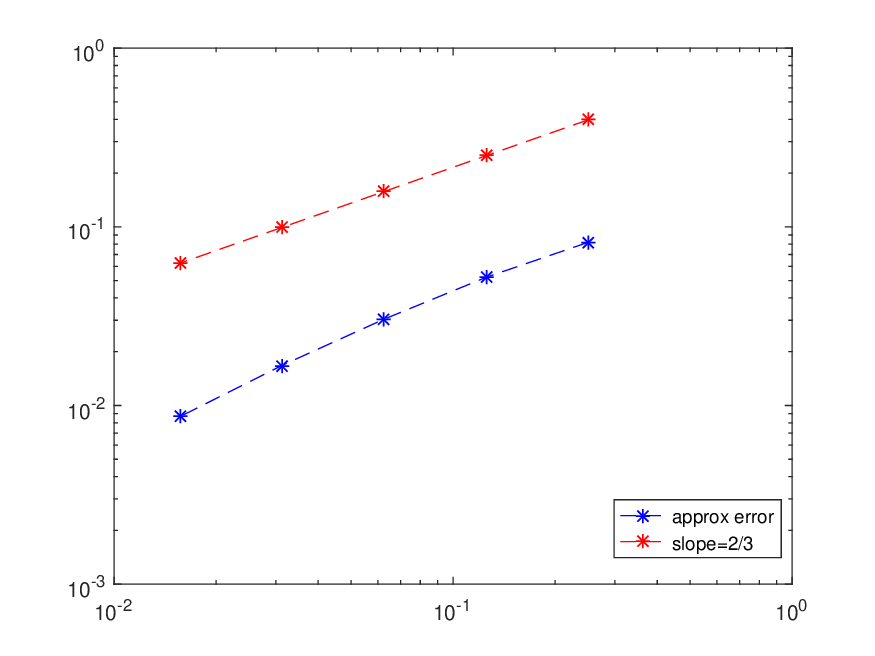}}
 \caption{Mean-square convergence rates for temporal discretization and space discretization}
\end{figure}

In view of Lemma \ref{eq-lem:discrete-trace-formula-full}, the energy evolution law $\mathbb{E}(J(U_h^n,V_h^n)$
grows linearly as time increases. 
In Figure 2, we show the evolution of discrete energy by  using the proposed methods \eqref{eq:U-formulation}-\eqref{eq:V-formulation} and one can detect the expected linear growth.  Here the expectation is approximated by taking average over 1000 realizations with  step-size $k=2^{-8}$ and $h=2^{-8}$.
\begin{figure}[!ht]
\centering
\scalebox{0.5}[0.51]{\includegraphics[0.5in,
0.2in][5.3in,4.1in]{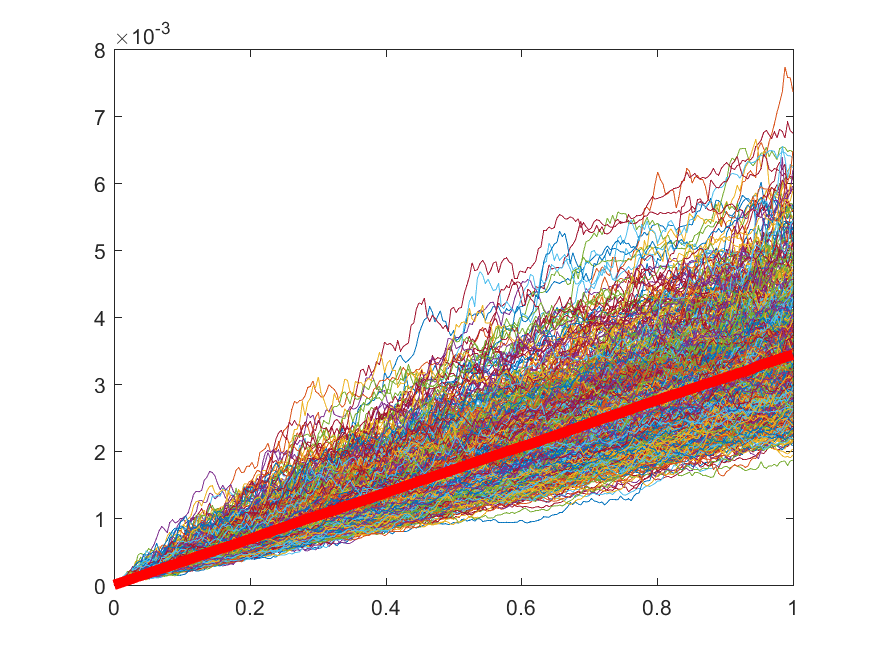}}
 \caption{Energy-preserving property }
\end{figure}
\section{Appendix: Proof of Lemma \ref{lem:eq-spatial-temporal-regularity}}

Recalling the definition of the Lyaponov function $J(\cdot,\cdot)$ in Theorem \ref{them:existen=expon-integ-result-mild}, we apply Ito formula to $J(u_h,v_h)$ to get
\begin{align}
\begin{split}
J(u_h(t),v_h(t))
=
&
J(u_h(0),v_h(0))
+
\int_0^t\big<DJ(u_h,v_h),\mathbb{F}_h(X_h)+A_hX_h\big>\,\dd s
\\
&
+
\frac12\int_0^t\mathrm{Tr}\left((\mathbb{B}_hQ^{\frac12})(\mathbb{B}_hQ^{\frac12})^*D^2J(u_h,v_h)\right)\,\dd s
+
\int_0^t\big<(DJ(u_h,v_h)),\mathbb{B}_h\dd W(s)\big>
\\
=&
J(u_h(0),v_h(0))
+
\frac12\int_0^t\mathrm{Tr}(P_hQP_h)\,\dd s
+
\int_0^t\big<v_h,P_h\,\dd W(s)\big>.
\end{split}
\end{align}
We take the pth's power, the  supremum with respect
to $s\in [0,T]$ and then the expectation and also use  the
Burkholder-Davis-Gundy inequality (\cite[Lemma 7.2]{da2014stochastic}), to obtain
\begin{align}
\begin{split}
&\mathbb{E}\sup_{s\in[0,T]}[J^p(u_h(s),v_h(s))]
\\
\leq
&
C\mathbb{E}[J^p(P_hu_0,P_hv_0)]
+
C\|P_hQ^\frac12\|_{\mathcal{L}_2(H)}^{2p}
+
C\mathbb{E}\sup_{s\in[0,T]}\Big[\Big\|\int_0^s\big<v_h(r),\,P_h\dd W(r)\big>\Big\|^p\Big]
\\\
\leq
&
C\mathbb{E}[J^p(P_hu_0,P_hv_0)]
+
C
\|Q^\frac12\|_{\mathcal{L}_2(H)}^{2p}
+
C\mathbb{E}\Big[\Big(\int_0^T\|Q^\frac12v_h\|^2\,\dd s\Big)^{\frac p2}\Big]
\\
\leq
&
C\mathbb{E}[J^p(P_hu_0,P_hv_0)]
+
C
\|Q^\frac12\|_{\mathcal{L}_2(H)}^{2p}
+
C\|Q^\frac12\|_{\mathcal{L}(H)}\int_0^T\mathbb{E}\sup_{r\in[0,s]}[J^p(u_h(r),v_h(r))]\dd s.
\end{split}
\end{align}
Hence the inequality \eqref{lem:eq-temporal-regularity-semi} follows from Gronwall's lemma.
\begin{align}
\mathbb{E}\sup_{s\in[0,T]}[J^p(u_h(s),v_h(s))]
<\infty.
\end{align}

To show \eqref{lem:eq-exponential-integ-properties-semi}, we define a linear operator as follows
\begin{align}
\begin{split}
(\mathbb{D}_{A_h,\mathbb{F}_h,\mathbb{B}_h}J)(u_h,v_h)):
=
&
\big<D_{u_h}J(u_h,v_h), v_h\big>_{L_2}+\big<D_{v_h}J(u_h,v_h), \Lambda_hu_h-P_hf(u_h)\big>_{L_2}
\\
&+
\frac12 \sum_{i=1}^\infty\big<D_{v_h,v_h}J(u_h,v_h)P_hQ^{\frac12}e_i,P_hQ^{\frac12}e_i\big>.
\end{split}
\end{align}
Then we arrive at
\begin{align}
\begin{split}
\mathbb{D}_{A_h,\mathbb{F}_h,\mathbb{B}_h}(X_h)+
\frac1{2\exp(\alpha t)}\|B_h^*(DJ(u_h,v_h))\|^2
\leq
&
\frac12\mathrm{Tr}(Q)
+
\frac1{2\exp(\alpha t)}\sum_{i=1}^\infty\big<v_h(s),Q^{\frac12}e_i\big>_{L^2}^2
\\
\leq
&
\frac12\mathrm{Tr}(Q)+
\frac1{2\exp(\alpha t)}\|v_h(s)\|^2\mathrm{Tr}(Q),
\end{split}
\end{align}
which together with exponential integrability lemma in  \cite[Corollary 2.4]{Cox2013Local} yields
\begin{align}\label{eq:bound-semi-discrete-solution-integr}
\mathbb{E}\Big(\exp\Big(\frac{J(u_h(t),v_h(t))}{\exp(\alpha t)}-\int_0^t\frac{\mathrm{Tr}Q}{2\exp(\alpha s)}\dd s\Big)\Big)
\leq
\mathbb{E}\Big(\exp(J(P_hu_0,P_hv_0)\Big)
<\infty.
\end{align}
Therefore, this in combination with  Jensen's inequality, Young's inequality and the fact $\|u_h\|_{L^6}\leq \|u\|^{\frac a2}\|A^{\frac12}u_h\|^{1-a}$ with $a=\frac d3$  implies
\begin{align}
\begin{split}
\mathbb{E}\Big[\exp\Big(\int_0^tc\|u_h(s)\|_{L^6}^2\,\dd s\Big)\Big]
\leq
&
C(T)\sup_{s\in[0,T]}\mathbb{E}\big(\exp(cT\|u_h(s)\|_{L^6}^2)\big)
\\
\leq
&
C(T)\sup_{s\in[0,T]}\mathbb{E}\Big(\exp\Big(\frac{\|A^{\frac12}u_h(s)\|^2}{2\exp{\alpha s}}\Big)
\exp\Big(\exp(\frac a{1-a} \alpha T  )\|u_h(s)\|^2(CT)^{\frac 1{1-a}}2^{\frac a{1-a}}\Big)\Big)
\\
\leq
&
C\sup_{s\in[0,T]}\mathbb{E}\Big(\exp\Big(\frac{\|\nabla u_h(s)\|^2}{2\exp(\alpha s)}\Big) \exp\Big(C(T)\|u_h(s)\|^2\Big)
\\
\leq
&
C\sup_{s\in[0,T]}\mathbb{E}\Big(\exp\Big(\frac{\|\nabla u_h(s)\|^2}{2\exp(\alpha s)}\Big) \exp\Big(\frac{\|u_h(s)\|^4}{2\exp(\alpha s)}+C(T)\Big)
\\
<&
\infty.
\end{split}
\end{align}
Hence this finishes the proof of this lemma.  $\hfill\square$
\bibliography{bibfile}

\bibliographystyle{abbrv}
 \end{document}